\documentclass[11pt]{article}

\usepackage{amssymb,epsfig}
\usepackage{amsmath}
\usepackage{comment}
\usepackage[dvips]{color}
\usepackage[T1]{fontenc}
\usepackage[latin1]{inputenc}
\usepackage[textwidth=2.5cm]{todonotes}

\newtheorem{defi}{Definition}
\newtheorem{prop}[defi]{Proposition}
\newtheorem{theo}[defi]{Theorem}
\newtheorem{conj}[defi]{Conjecture}
\newtheorem{lemm}[defi]{Lemma}
\newtheorem{coro}[defi]{Corollary}
\newtheorem{rema}[defi]{Remark}
\newtheorem{exem}[defi]{Example}
\newtheorem{exems}[defi]{Examples}

\newcommand{\bdefi}{\begin{defi}}
\newcommand{\edefi}{\end{defi}}
\newcommand{\bprop}{\begin{prop}}
\newcommand{\eprop}{\end{prop}}
\newcommand{\btheo}{\begin{theo}}
\newcommand{\etheo}{\end{theo}}
\newcommand{\blemm}{\begin{lemm}}
\newcommand{\brema}{\begin{rema}}
\newcommand{\erema}{\end{rema}}
\newcommand{\bexer}{\begin{exem}}
\newcommand{\eexer}{\end{exem}}
\newcommand{\bexems}{\begin{exems}}
\newcommand{\eexems}{\end{exems}}
\newcommand{\bconj}{\begin{conj}}
\newcommand{\econj}{\end{conj}}
\newcommand{\elemm}{\end{lemm}}
\newcommand{\bcoro}{\begin{coro}}
\newcommand{\ecoro}{\end{coro}}
\newcommand{\dem}{\noindent{\bf Proof. }}

\usepackage{mathrsfs}
\renewcommand\mathcal{\mathscr}

\newcommand{\N}{{\cal N}}

\newcommand{\D}{{\cal D}}

\newcommand{\F}{{\cal F}}
\newcommand{\V}{{\cal V}}
\newcommand{\W}{{\cal W}}

\renewcommand{\H}{{\cal H}}

\newcommand{\C}{{\cal C}}

\newcommand{\maths}[1]{{\mathbb #1}}  

\newcommand{\RR}{\maths{R}}
\newcommand{\NN}{\maths{N}}

\newcommand{\HH}{\maths{H}}

\newcommand{\ZZ}{\maths{Z}}

\newcommand{\weakstar}{\overset{*}\rightharpoonup}
\newcommand{\ra}{\rightarrow}
\newcommand{\bs}{\backslash}

\newcommand{\ov}[1]{{\overline #1}} 
\newcommand{\wt}[1]{{\widetilde{#1}}}

\newcommand{\ga}{\gamma}
\newcommand{\Ga}{\Gamma}

\newcommand{\cqfd}{\hfill$\Box$}

\newcommand{\card}{{\operatorname{Card}}}

\newcommand{\Par}{\operatorname{Par}}

\newcommand{\CAT}{\operatorname{CAT}}

\newcommand{\stab}{\operatorname{Stab}}
\newcommand{\flow}[1]{{g^{#1}}}  
\newcommand{\muss}[1]{{\mu^{\rm ss}_{#1}}}  
\newcommand{\mus}[1]{{\mu^{\rm s}_{#1}}}  
\newcommand{\musu}[1]{{\mu^{\rm su}_{#1}}}
\newcommand{\muu}[1]{{\mu^{\rm u}_{#1}}}  
\newcommand{\normal}[1]{\partial^1_{+}{#1}}

\newcounter{fig}

\usepackage{enumitem}
\setlist{nolistsep}

\title{Skinning measures in negative curvature \\
and equidistribution of equidistant submanifolds}
\author{Jouni Parkkonen \and Fr\'ed\'eric Paulin} 
\date{}

\begin{document}
\bibliographystyle{../alphanum}

\maketitle

\begin{abstract} Let $C$ be a locally convex subset of a negatively
  curved Riemannian manifold $M$. We define the skinning measure
  $\sigma_C$ on the outer unit normal bundle to $C$ in $M$ by pulling
  back the Patterson-Sullivan measures at infinity, and give a
  finiteness result of $\sigma_C$, generalising the work of Oh and
  Shah, with different methods. We prove that the skinning measures,
  when finite, of the equidistant hypersurfaces to $C$ equidistribute
  to the Bowen-Margulis measure $m_{\rm BM}$ on $T^1M$, assuming only
  $m_{\rm BM}$ is finite and mixing for the geodesic flow. Under
  additional assumptions on the rate of mixing, we give a control on
  the rate of equidistribution.  
\end{abstract}

\section{Introduction}

Let $M$ be a complete connected Riemannian manifold with
sectional curvature at most $-1$.  For any proper nonempty properly
immersed locally convex subset $C$ of $M$ and $t>0$, let $\Sigma_{t}$
be the (Lipschitz) submanifold of $T^1M$ that consists of images by
the geodesic flow at time $t$ of the outward-pointing unit normal
vectors to the boundary of $C$ (see Section \ref{sec:geometry} for
precise definitions).

If $M$ has constant curvature and finite volume and if $C$ is an
immersed totally geodesic submanifold of finite volume,  we showed
in \cite[Thm.~2.2.]{ParPau12JMD} that the Riemannian measure of
$\Sigma_{t}$ equidistributes to the Liouville measure of $T^1M$ (which
is the Riemannian measure of the Sasaki metric of $T^1M$). This result
also follows from the equidistribution result of Eskin and McMullen
\cite[Thm.~1.2]{EskMcMul93} in affine symmetric spaces, see \cite[\S
4]{ParPauRev} for details.

In this paper, we generalise the above result when $C$ is no longer
required to be totally geodesic and when $M$ has variable
curvature. Though the methods of locally homogeneous spaces as in
\cite{EskMcMul93} are then completely not applicable, the strategy of
\cite{ParPau12JMD} remains helpful. Both the measures on $T^1M$ and on
$\Sigma_{t}$ need to be adapted to variable curvature.

The measure on $T^1M$ we will consider (when $M$ is nonelementary and
its fundamental group has finite critical exponent) is the well-known
{\em Bowen-Margulis measure} $m_{\rm BM}$ (see \cite{Roblin03} for a
nice presentation). It coincides with the Liouville measure (up to a
multiplicative constant) when $M$ is locally symmetric with finite
volume (see for instance \cite[\S 7]{ParPauRev} when $M$ is real
hyperbolic). It is, when finite and normalised, the unique probability
measure of maximal entropy for the geodesic flow on $T^1M$ (see
\cite{Margulis70} and \cite{Bowen72a} when $M$ is compact, and
\cite{OtaPei04} under the only assumption that $m_{\rm BM}$ is
finite). The Bowen-Margulis measure is finite for instance when $M$ is
compact, or when $M$ is geometrically finite and the critical exponent
of its fundamental group is strictly bigger than the critical
exponents of its parabolic subgroups (as it is the case when $M$ is
locally symmetric), by \cite{DalOtaPei00}.  By
\cite[Thm.~1]{Babillot02b}, the Bowen-Margulis measure, when finite,
is mixing if the length spectrum of $M$ is not contained in a discrete
subgroup of $\RR$. By \cite{Dalbo99,Dalbo00}, this condition holds for
instance when $M$ is $2$-dimensional or locally symmetric, or if its
fundamental group contains a parabolic element.

The measure on $\Sigma_{t}$ we will consider is the {\em skinning
  measure} that we introduce under this generality in this paper (see
Section \ref{sec:skinning}), as an appropriate pushforward to
$\Sigma_{t}$ of the natural measures at infinity of the universal
cover of $M$. It scales by $e^{\delta t}$, where $\delta$ is the
critical exponent of $M$, under the geodesic flow map from
$\Sigma_{s}$ to $\Sigma_{s+t}$.  When $C$ is an immersed horoball,
$\Sigma_{t}$ is a leaf of the strong unstable foliation of the
geodesic flow on $T^1M$, and the skinning measure on $\Sigma_{t}$ is
simply the conditional measure of the Bowen-Margulis measure on this
leaf (see for instance \cite{Margulis04,Roblin03}).  When $M$ is
geometrically finite with constant curvature, and when $C$ is an
immersed ball, horoball or totally geodesic submanifold, the skinning
measure on $\Sigma_{t}$ has been introduced by Oh and Shah
\cite{OhShaInv,OhShaCounting}, who coined the term, with beautiful
applications to circle packings, and coincides with the Riemannian
measure up to a multiplicative constant (see \cite[\S 7]{ParPauRev}
for a computation of the constant) when furthermore $M$ has finite
volume.  When the intersection of $\Sigma_{t}$ with the nonwandering
set of the geodesic flow of $T^1M$ is compact, the skinning measure is
finite.

When $M$ is geometrically finite, generalizing (and giving an
alternative proof of) Theorem 6.4 in \cite{OhShaCounting} which
assumes the curvature to be constant, we give in Theorem
\ref{theo:geneohshah} a sharp criterion for the finiteness of the
skinning measure, by studing its decay in the cusps of $M$. This decay
is analogous to the decay of the Bowen-Margulis measure in the cusps,
which was first studied by Sullivan \cite{Sullivan84} who called it
the fluctuating density property (see also \cite{StrVel95} and
\cite[Theo.~4.1]{HerPau04}). The criterion, as in the case of the
Bowen-Margulis measure in \cite{DalOtaPei00}, is a separation property
of the critical exponents.

The following theorem is a simplified version of the main result of
this paper. In the more general result, Theorem \ref{theo:equid} in
Section \ref{sec:equidistribution}, we replace $\Sigma_t$ by
$\flow{t}\Omega$, where $\Omega$ is an open set of outward-pointing
unit normal vectors to $\partial C$ with finite nonzero skinning
measure.

\btheo\label{theo:intro} Let $M$ be a nonelementary connected complete
Riemannian manifold with pinched negative sectional curvature.  Assume
that the Bowen-Margulis measure on $T^1M$ is finite and mixing for the
geodesic flow.  Let $C$ be a proper nonempty properly immersed locally
convex subset of $M$ with finite nonzero skinning measure. Then as $t$
tends to $+\infty$, the skinning measure on $\Sigma_{t}$
equidistributes to the Bowen-Margulis measure on $T^1M$.  
\etheo

When $C$ is an immersed ball or horoball, then this result is due to
Margulis when $M$ has finite volume, see for example \cite{Margulis04},
and to Babillot \cite[Theo.~3]{Babillot02b} and Roblin \cite{Roblin03}
under the weak assumptions of Theorem $1$. Many ideas of our proof go
back to \cite{Margulis69}. See also \cite{Schapira04b, Marklof10,
  KonOh11, Kim13} for other results on the equidistribution of
horospheres and applications.

For instance, it follows from Theorem \ref{theo:intro} that when $M$
is a compact Riemannian manifold with negative sectional curvature,
when $C$ is the image in $M$ of the convex hull of the limit set of a
convex-cocompact subgroup of the covering group of a universal cover
of $M$, then the skinning measure on $\Sigma_{t}$ equidistributes to
the Bowen-Margulis measure on $T^1M$. But we make no compactness
assumption in our theorem, only requiring the finiteness of the
measures under consideration. The main tool is a general
disintegration result of the Bowen-Margulis measure over any skinning
measure (see Proposition \ref{prop:disintegration}).

\medskip We also give (see Section \ref{sec:expo}) estimates on the
rate of equidistribution in the previous result, under assumptions on
the rate of mixing of the geodesic flow. When $M$ is locally symmetric
and arithmetic, the rate of mixing of the geodesic flow for
sufficiently smooth functions is exponential, by the work of Kleinbock
and Margulis \cite[Theo.~2.4.3]{KleMar96} and Clozel
\cite[Theo.~3.1]{Clozel03}. When the curvature is variable, the
appropriate regularity is the H\"older one. The rate of mixing of the
geodesic flow for H\"older-continuous functions is exponential if $M$
is compact and has dimension $2$ by the work of Dolgopyat
\cite{Dolgopyat98} or if $M$ is compact and locally symmetric (without
the arithmetic assumption) by \cite[Coro.~1.5]{Stoyanov11} (see also
\cite{Liverani04} when $M$ is compact, the result stated for
the Liouville measure should extend to the Bowen-Margulis measure, for
instance by using the tools of \cite{GiuLivPol12} if the sectional
curvature of $M$ is $\frac{1}{9}$-pinched).

\btheo\label{theo:introexpo} Under the hypotheses of Theorem
\ref{theo:intro}, in anyone of the above cases when the geodesic flow
of $T^1M$ is mixing with exponential speed, the skinning measure
$\sigma_t$ of $\Sigma_{t}$ equidistributes to the Bowen-Margulis
measure with exponential speed.  \etheo

More precisely in the H\"older-case, if $M$ is compact and is
$2$-dimensional or locally symmetric, then there exist $\alpha\in
\;]0,1[\,$ and $\tau>0$ such that for every
$\alpha$-H\"older-continuous function $\psi:T^1M\ra\RR$
with $\alpha$-H\"older norm $\|\psi\|_\alpha$ (see Section
\ref{sec:expo} for precise definitions), as $t$ tends to $+\infty$,
$$
\frac{1}{\|\sigma_t\|}\int_{\Sigma_t} \psi \;d\sigma_t=
\frac{1}{\|m_{\rm BM}\|}\int_{T^1M}
\psi\;dm_{\rm BM}+O(e^{-\tau\, t}\;\|\psi\|_\alpha)\;.
$$

In \cite{ParPau13b}, we will use the tools introduced in this paper to
study counting results of common perpendicular arcs between locally
convex subsets in variable negative curvature.

\medskip\noindent{\small {\it Acknowledgement: } The first author
  thanks the University of Paris-Sud (Orsay) for a month of visiting
  professor where this work was started, and the FIM of ETH Z\"urich
  for its support in 2011-2012 when this work was completed.  The
  second author thanks the ETH in Z\"urich for frequent secret stays
  during the completion of the writing of this paper.  We thank
  P.~Koskela, K.~Rajala and H.~Tuominen for useful analytical
  conversations.  We thank the referee for numerous very helpful
  remarks, that have greatly improved this paper.}

\section{Geometry, dynamics and convexity in negative 
curvature}
\label{sec:geometry}

In this section, we review briefly the required background on
negatively curved Riemannian manifolds, seen as locally
$\CAT(-\kappa)$ spaces, using for instance \cite{BriHae99} as a
general reference, and their unit tangent bundles and geodesic flows.
We introduce the geometric fibred neighbourhoods of the outer unit
normal bundle of the boundary of a convex subset that will be used in
what follows.

\medskip
\noindent{\bf Geometry and dynamics.} Let $\wt M$ be a complete simply
connected Riemannian manifold with 
sectional curvature bounded above by $-1$, and let $x_{0}\in\wt M$.
Let $\Ga$ be a discrete, nonelementary group of isometries of $\wt M$,
and let us denote the quotient space of $\wt M$ under $\Ga$ by
$M=\Ga\bs \wt M$.  We denote by $\partial_{\infty}\wt M$ the boundary
at infinity of $\wt M$ (with its usual H\"older structure), by
$\Lambda\Ga$ the limit set of $\Ga$ and by $\C\Lambda\Ga$ the convex
hull in $\wt M$ of $\Lambda\Ga$. For every $\epsilon>0$, we denote by
$\N_\epsilon A$ the closed $\epsilon$-neighbourhood of a subset $A$ of
$\wt M$, and by convention $\N_0A=\overline{A}$.

For any point $\xi\in\partial_{\infty}\wt M$, let
$\rho_{\xi}:[0,+\infty[\ra \wt M$ be the geodesic ray with origin
$x_{0}$ and point at infinity $\xi$.  The {\em Busemann cocycle} of
$\wt M$ is the map $\beta: \wt M\times\wt
M\times\partial_{\infty} \wt M\to\RR$ defined by
$$
(x,y,\xi)\mapsto \beta_{\xi}(x,y)=
\lim_{t\to+\infty}d(\rho_{\xi}(t),x)-d(\rho_{\xi}(t),y)\;.
$$
The above limit exists and is independent of $x_{0}$. If $y$ is a
point in the (image of the) geodesic ray from $x$ to $\xi$, then
$\beta_{\xi}(x,y)=d(x,y)$.  The Busemann cocycle satisfies
\begin{equation}\label{eq:cocycle}
\beta_{\ga \xi}(\ga x,\ga y)=\beta_{\xi}(x,y)\;\;\;{\rm and}\;\;\;
\beta_{\xi}(x,y)+\beta_{\xi}(y,z)=\beta_{\xi}(x,z)\;,
\end{equation}
for all $\xi\in\partial_{\infty}\wt M$, all $x,y,z\in\wt M$ and every
isometry $\ga$ of $\wt M$.  The {\em visual distance} $d_{x_0}$ (based
at $x_0$) on $\partial_{\infty}\wt M$ is the distance defined by 
\begin{equation}\label{eq:visdist}
d_{x_{0}}(\xi,\eta)=
e^{-\frac 12(\beta_{\xi}(x_{0},y)+\beta_{\eta}(x_{0},y))}
\end{equation}
for any $y$ in the geodesic line between $\xi $ and $\eta$ if
$\xi\ne\eta$, and $d_{x_{0}}(\xi,\eta)=0$ if $\xi=\eta$.

The unit tangent bundle $T^1N$ of a complete Riemannian manifold $N$
can be identified with the set of locally geodesic lines $\ell:\RR\to
N$ in $N$, endowed with the compact-open topology. More precisely, we
identify a locally geodesic line $\ell$ and its (unit) tangent vector
$\dot{\ell}(0)$ at time $t=0$ and, conversely, any $v\in T^1 N$ is the
tangent vector at time $t=0$ of a unique locally geodesic line.  We
will use this identification without mention in this paper. We denote
by $\pi: T^1 N\to N$ the base point projection, which is given by
$\pi(\ell)=\ell(0)$.

The {\em geodesic flow} on $T^1N$ is the dynamical system
$(\flow{t})_{t\in\RR}$, where $\flow{t}\ell\,(s)=\ell(s+t)$, for all
$\ell\in T^1N$ and $s,t\in\RR$.  The isometry group of $\wt M$ acts on
the space of geodesic lines in $\wt M$ by postcomposition: $(\ga,\ell)
\mapsto \ga\circ\ell$, and this action commutes with the geodesic
flow.

When $\Ga$ acts on $\wt M$ without fixed point, we have an
identification $\Ga\bs T^1\wt M=T^1M$.  Even in the general case with
torsion, we denote by $T^1M$ the quotient space $\Ga\bs T^1\wt M$. We
use the notation $(\flow{t})_{t\in\RR}$ also for the geodesic flow on
$T^1 M$ (induced by the geodesic flow on $T^1\wt M$ by passing to the
quotient).

We denote by $\iota:T^1\wt M\ra T^1\wt M$ the {\it antipodal (flip)
  map} $v\mapsto -v$, and we again denote by $\iota:T^1M\ra T^1M$ its
quotient map. We have $\iota\circ \flow t= \flow{-t}\circ \iota$.

For every unit tangent vector $v\in T^1\wt M$, let $v_{-}=v(-\infty)$
and $v_{+}=v(+\infty)$ be the two endpoints in the sphere at infinity
of the geodesic line defined by $v$. Let $\partial_{\infty}^2\wt M$ be
the subset of $\partial_{\infty}\wt M\times\partial_{\infty}\wt M $
which consists of pairs of distinct points at infinity.  The {\em Hopf
  parametrisation} of $T^1\wt M$ is the identification of $v\in T^1\wt
M$ with the triple $(v_{-},v_{+},t) \in \partial_{\infty}^2\wt M\times
\RR$, where $t$ is the signed (algebraic) distance of $\pi(v)$ from
the closest point $p_{v,x_0}$ to $x_{0}$ on the (oriented) geodesic
line defined by $v$. This map is a homeomorphism, the geodesic flow
acts by $\flow s(v_{-},v_{+},t) =(v_{-},v_{+},t+s)$ and, for every
isometry $\ga$ of $\wt M$, the image of $\ga v$ is $(\ga v_{-}, \ga
v_{+}, t+ t_{\ga,v_-,v_+})$, where $t_{\ga,v_-,v_+}$ is the signed
distance from $\ga p_{v,x_0}$ to $p_{\ga v,x_0}$.  Furthermore, in
these coordinates, the antipodal map $\iota$ is $(v_-,v_+,t)\mapsto
(v_+,v_-,-t)$.

The {\em strong stable manifold} of $v\in T^1\wt M$ is 
$$
W^{\rm ss}(v)=\{v'\in T^1\wt M:d(v(t),v'(t))\to
0 \textrm{ as } t\to+\infty\},  
$$
and  its {\em strong unstable manifold}  is 
$$
W^{\rm su}(v)=\{v'\in T^1\wt M:d(v(t),v'(t))\to 0 
\textrm{ as } t\to-\infty\}, 
$$
The union for $t\in\RR$ of the images under $\flow t$ of the strong
stable manifold of $v\in T^{1}\wt M$ is the {\em stable manifold}
$W^{\rm s}(v)=\bigcup_{t\in\RR}\flow t W^{\rm ss}(v)$ of $v$, which
consists of the elements $v'\in T^1\wt M$ with $v'_+=v_+$. Similarly,
the union of the images under the geodesic flow at all times of the
strong unstable manifold of $v$ is the {\em unstable manifold} $W^{\rm
  u}(v)$ of $v$, which consists of the elements $v'\in T^1\wt M$ with
$v'_-=v_-$.

The strong stable manifolds, stable manifolds, strong unstable
manifolds and unstable manifolds are the (smooth) leaves of
foliations, that are invariant under the geodesic flow and the
isometry group of $\wt M$, denoted by $\W^{\rm ss}, \W^{\rm s},
\W^{\rm su}$ and $\W^{\rm u}$, respectively.  These foliations are
H\"older-continuous when $\wt M$ has pinched negative sectional
curvature with bounded derivatives (see for instance \cite{Brin95},
\cite[\S 7.1]{PauPolSha11}).  The maps from $\RR\times W^{\rm ss}(v)$
to $W^{\rm s}(v)$ defined by $(s,v')\mapsto \flow sv'$ and from
$\RR\times W^{\rm su}(v)$ to $W^{\rm u}(v)$ defined by $(s,v')\mapsto
\flow sv'$ are smooth diffeomorphisms.

\smallskip\noindent
\begin{minipage}{9.9cm} ~~~ The images of the strong unstable and
  strong stable manifolds of $v\in T^1\wt M$ under the base point
  projection, denoted by $H_{-}(v)= \pi(W^{\rm su}(v))$ and $H_{+}(v)=
  \pi(W^{\rm ss}(v))$, are called, respectively, the {\em unstable and
    stable horospheres} of $v$, and are said to be {\it centered at}
  $v_{-}$ and $v_{+}$, respectively. The unstable horosphere of $v$
  coincides with the zero set of the map $x\mapsto f_{-}(x)=
  \beta_{v_-}(x,\pi(v))$, and, similarly, the stable horosphere of $v$
  coincides with the zero set of $x\mapsto f_{+}(x)= \beta_{v_+}(x,
  \pi(v))$.  The corresponding sublevel sets $H\!B_{-}(v)=
  f_{-}^{-1}(]-\infty,0])$ and $H\!B_{+}(v)= f_{+}^{-1}(]-\infty,0])$
  are called the {\em horoballs} bounded by $H_{-}(v)$ and $H_{+}(v)$.
  Horoballs are (strictly) convex subsets of $\wt M$.
\end{minipage}
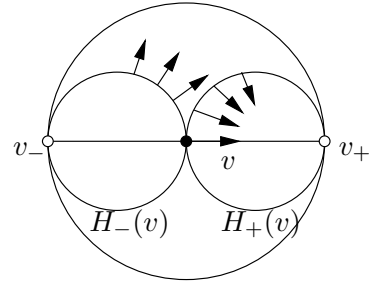
\begin{minipage}{5cm}
\begin{center}
\input{fig_horobull.pstex_t}
\end{center}
\end{minipage}

\smallskip For every $w\in T^1\wt M$, let $d_{W^{\rm ss}(w)}$ be the {\it
  Hamenst\"adt distance} on the strong stable leaf of $w$, defined
as follows (see \cite{Hamenstadt89}, \cite[Appendix]{HerPau97}, as well
as \cite[\S 2.2]{HerPau10} for a generalisation when the horosphere
$H_{+}(w)$ is replaced by the boundary of any nonempty closed convex
subset): for all $v,v'\in W^{\rm ss}(w)$,
$$
d_{W^{\rm ss}(w)}(v,v') = 
\lim_{t\ra+\infty} e^{\frac{1}{2}d(v(-t),\;v'(-t))-t}\;.
$$
This limit exists, and the Hamenst\"adt distance is a distance inducing
the original topology on $W^{\rm ss}(w)$. For all $v,v'\in W^{\rm ss}
(w)$ and for every isometry $\ga$ of $\wt M$, we have $d_{W^{\rm ss}
  (\ga w)}(\ga v,\ga v')= d_{W^{\rm ss}(w)}(v,v')$.  By the triangle
inequality, for all $v,v'\in W^{\rm ss} (w)$, we have
\begin{equation}\label{eq:majodistHamen}
d_{W^{\rm ss} (w)}(v,v')\leq e^{\frac{1}{2}d(\pi(v),\,\pi(v'))}\;.
\end{equation}
For all $w\in
T^1\wt M$, $s\in\RR$ and $v,v'\in W^{\rm ss} (w)$, we have
\begin{equation}\label{eq:contracHamdist}
d_{W^{\rm ss} (\flow sw)}(\flow sv,\flow sv')=e^{-s}d_{W^{\rm ss}(w)}(v,v')\;.
\end{equation}

\medskip A usual distance $d$ on $T^1\wt M$ is defined, 
for all $v,v'\in T^1\wt M$, by
$$
d(v,v')=\frac{1}{\sqrt{\pi}}\;\int_\RR d(v(t),v'(t))\,e^{-t^2}\;dt\;.
$$
This distance is invariant under the group of isometries of $\wt M$
and the antipodal map. Also note that for all $s\in\RR$ and $v\in
T^1\wt M$, we have
\begin{equation}\label{eq:distlelongeod}
d(g^sv,v)=s\;.
\end{equation}

\blemm \label{lem:compardisttunwss} There exists $c>0$ such that for
all $w\in T^1\wt M$ and $v,v'\in W^{\rm ss}(w)$, we have
$$
d(v,v')\leq c\;d_{W^{\rm ss}(w)}(v,v')\;.
$$
\elemm

\dem We may assume that $v\neq v'$. By the convexity properties of the
distance in $\wt M$, the map from $\RR$ to $\RR$ defined by $t\mapsto
d(v(t),v'(t))$ is decreasing, with image $]0,+\infty[$. Let $S\in\RR$
be such that $d(v(S),v'(S))=1$.  For every $t\le S$, let $p$ and $p'$
be the closest point projections of $v(S)$ and $v'(S)$ on the geodesic
segment $[v(t),v'(t)]$. We have $d(p,v(S)),d(p',v'(S))\leq 1$ by
comparison. Hence, by convexity and the triangle inequality,
\begin{align*}
d(v(t),v'(t))&\geq d(v(t),p)+d(p',v'(t))\\& 
\geq d(v(t),v(S))-1+ d(v'(t),v'(S))-1=2(S-t-1)\;.
\end{align*}
Thus by the definition of the Hamenst\"adt distance $d_{W^{\rm ss}(w)}$,
we have
\begin{equation}\label{eq:minodistss}
d_{W^{\rm ss}(w)}(v,v')\geq e^{S-1}\;.
\end{equation}

By the triangle inequality, if $t\leq S$, then
$$
d(v(t),v'(t))\leq
d(v(t),v(S))+d(v(S),v'(S))+ d(v'(S)),v'(t)=2(S-t)+1\;.
$$
Since $\wt M$ is $\CAT(-1)$, if $t\geq S$, we have by comparison
$$
d(v(t),v'(t))\leq e^{S-t}\;d(v(S),v'(S))=e^{S-t}\;.
$$
Therefore, by the definition of the distance $d$ on $T^1\wt M$,
$$
d(v,v')\leq \int_{-\infty}^{S}(2(S-t)+1)\,e^{-t^2}\;dt+
\int_{S}^{+\infty}e^{S-t}\,e^{-t^2}\;dt=\operatorname{O}(e^{S})\;.
$$
The result hence follows from Equation \eqref{eq:minodistss}.  
\cqfd

\bigskip \noindent{\bf Convexity.} Let $C$ be a nonempty closed convex
subset of $\wt M$.  We denote by $\partial C$ the boundary of $C$ in
$\wt M$ and by $\partial_{\infty}C$ its set of points at infinity (the
set of endpoints of geodesic rays contained in $C$). Let $P_{C}:\wt
M\cup (\partial_{\infty} \wt M -\partial_{\infty}C)\to C$ be the
(continuous) closest point map: if $\xi\in\partial_{\infty}\wt
M-\partial_{\infty}C$, then $P_{C}(\xi)$ is defined to be the unique
point in $C$ that minimises the map $x\mapsto\beta_{\xi}(x,x_{0})$
from $C$ to $\RR$. For every isometry $\ga$ of $\wt M$, we have
$P_{\ga C}\circ\ga=\ga\circ P_C$.

Let $\normal C$ be the subset of $T^1\wt M$ consisting of the geodesic
lines $v:\RR\to \wt M$ with $v(0)\in\partial C$, $v_+
\notin \partial_\infty C$ and $P_{C}(v_+) =v(0)$. Note that
$\pi(\normal C)=\partial C$ and that for every isometry $\ga$ of $\wt
M$, we have $\normal (\ga C)=\ga\,\normal C$. In particular, $\normal
C$ is invariant under the isometries of $\wt M$ that preserve
$C$. When $C=H\!B_-(v)$ is the unstable horoball of $v\in T^1\wt M$,
then $\normal C$ is the strong unstable manifold $W^{\rm su}(v)$ of
$v$, and similarly, $W^{\rm ss}(v)=\iota \,\normal H\!B_+(v)$.

The restriction of $P_{C}$ to $\partial_{\infty} \wt M
-\partial_{\infty} C$ (which is not necessarily injective) has a
natural lift to a homeomorphism
$$
\nu P_{C}:\partial_{\infty}\wt M-\partial_{\infty}C\to\normal C
$$
such that $\pi\circ\nu P_{C}=P_{C}$. The inverse of $\nu P_{C}$ is the
endpoint map $v\mapsto v_+$ from $\normal C$ to $\partial_{\infty} \wt
M - \partial_{\infty}C$. In particular, $\normal C$ is a topological
submanifold of $T^1\wt M$. For every $s\geq 0$, the geodesic flow
induces a homeomorphism $\flow s:\normal C\to\normal\N_{s}C$. For
every isometry $\ga$ of $\wt M$, we have $\nu P_{\ga C} \circ\ga=
\ga\circ \nu P_{C}$.  We refer for instance to \cite{Walter76} for the
notion of $\rm C^{1,1}$ and Lipschitz manifolds. When $C$ has nonempty
interior and $\rm C^{1,1}$ boundary, then $\normal C$ is the Lipschitz
submanifold of $T^1\wt M$ consisting of the outward-pointing unit
normal vectors to $\partial C$, and the map $P_{C}$ itself is a
homeomorphism (between $\partial_{\infty}\wt M-\partial_{\infty}C$ and
$\partial C$).  This holds, for instance, by \cite{Walter76}, when $C$
is the closed $\eta$-neighbourhood of any nonempty convex subset of
$\wt M$ with $\eta>0$.

We now define a canonical fundamental system of neighbourhoods, of
dynamical origin, of these outer unit normal bundles of boundaries of
convex sets. 
Let 
\begin{equation}\label{eq:defiUC}
U_C=\{v\in T^1\wt M:\ v_+\notin\partial_\infty C\}\;.
\end{equation}
Note that $U_C$ is an open subset of $T^1\wt M$, invariant under the
geodesic flow, which is empty if and only if $C=\wt M$, and is dense
in $T^1\wt M$ if the interior of $\partial_\infty C$ in
$\partial_\infty \wt M$ is empty. We have $U_{\ga C}=\ga U_C$ for
every isometry $\ga$ of $\wt M$ and, in particular, $U_C$ is invariant
under the isometries of $\wt M$ preserving $C$.

\medskip \noindent
\begin{minipage}{9.9cm} ~~~ Define a map $f_C:U_C\ra \normal C$, as
  the composition of the map from $U_C$ onto $ \partial_\infty \wt
  M-\partial_\infty C$ sending $v$ to $v_+$ and the homeomorphism $\nu
  P_{C}$ from $\partial_\infty \wt M-\partial_\infty C$ to $\normal
  C$. The map $f_C$ is a fibration as the composition of such a map
  with the homeomorphism $\nu P_{C}$. The fiber of $w\in \normal C$ is
  exactly its stable leaf $W^{\rm s}(w)=\{v\in T^1\wt
  M\;:\;v_+=w_+\}$. In particular, $U_C$ is the disjoint union of the
  leaves $W^{\rm s}(w)$ for $w\in \normal C$.
\end{minipage}
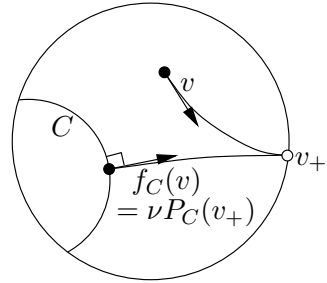
\begin{minipage}{5cm}
\begin{center}
\input{fig_UCfCnuPC.pstex_t}
\end{center}
\end{minipage}

\smallskip For every isometry $\ga$ of $\wt M$, we have $f_{\ga C}
\circ\ga=\ga\circ f_C$. We have $f_{\N_tC}= \flow t\circ f_C$ for all
$t\geq 0$, and $f_C\circ \flow t=f_C$ for all $t\in\RR$. In
particular, the fibration $f_C$ is invariant under the geodesic flow.

Let $\eta,R>0$. For all $w\in T^1M$, let
\begin{equation}\label{eq:defiboulehamen}
V_{w,\,R}=\{v'\in W^{\rm ss}(w)\;:\;d_{W^{\rm ss}(w)}(v',w)<R\}
\end{equation}
be the open ball of radius $R$ centered at $w$ for the Hamenst\"adt
distance in the strong stable leaf of $w$, and
\begin{align*}
V_{w,\,\eta,\,R}&=\{v\in W^{\rm s}(w)\;:\;\exists \;v'\in V_{w,\,R},
\;\exists\; s\in\;]-\eta,\eta\,[\;,\;\;\;\flow s v'=v\}\\ &=
\bigcup_{s\in\;]-\eta,\,\eta\,[}\flow s V_{w,\,R}=
\bigcup_{s\in\;]-\eta,\,\eta\,[}V_{\flow sw,\,e^{-s}R}\;.
\end{align*}
This last equality follows from the fact that, by Equation
\eqref{eq:contracHamdist}, we have $\flow sV_{w,\,R}= V_{\flow
  sw,\,e^{-s}R}$ for every $s\in\RR$. For every isometry $\ga$ of $\wt
M$, we have $\ga V_{w,\,R}= V_{\ga w,\,R}$ and $\ga V_{w,\,\eta,\,R}=
V_{\ga w,\,\eta,\,R}$. The map from $]-\eta,\eta[\;\times V_{w,\,R}$
to $V_{w,\,\eta,\,R}$ defined by $(s,v')\mapsto \flow s v'$ is a
homeomorphism.

\medskip \noindent
\begin{minipage}{9.9cm} ~~~ For every subset $\Omega$ of $\normal C$,
  let 
$$
\V_{\eta,\,R}(\Omega)=\bigcup_{w\in\Omega} V_{w,\,\eta,\,R}
$$
which is an open neighbourhood of $\Omega$ in $T^1\wt M$ if $\Omega$
is open in $\normal C$.  For every isometry $\ga$ of $\wt M$ and every
$t\geq 0$, we have $\ga\V_{\eta,\,R}(\Omega)=\V_{\eta,\,R}(\ga
\Omega)$ and
$$
\flow t\V_{\eta,\,R}(\Omega)=\V_{\eta,\,e^{-t}R}(\flow{t} \Omega)\;.
$$
\end{minipage}
\begin{minipage}{5cm}
\begin{center}
\input{fig_voisdyna.pstex_t}
\end{center}
\end{minipage}

\medskip These thickenings $\V_{\eta,\,R}(\Omega)$ are
nondecreasing in $\eta$ and in $R$ and  their intersection is $\Omega$.
Furthermore, we have 
$$\bigcup_{\substack{\eta>0\\R>0}}\V_{\eta,
  R}(\normal C)=U_C\,.$$ 
The restriction of $f_C$ to
$\V_{\eta,\,R}(\Omega)$ is a fibration over $\Omega$, with fiber of
$w\in \Omega$ the open subset $V_{w,\,\eta,\,R}$ of the stable leaf
of $w$.

\section{Patterson, Bowen-Margulis and skinning 
measures}
\label{sec:skinning}

Let $\wt M,\Ga,x_0, M$ and $T^1M$ be as in the beginning of Section
\ref{sec:geometry}. In this section, we first review some background
material on the measures associated with negatively curved manifolds
(for which we refer to \cite{Roblin03}). We then define the skinning
measure associated to any nonempty closed convex subset, generalising
the construction of \cite{OhShaCounting,OhShaInv}, and we prove some
basic properties of these measures, as well as a crucial
disintegration result. Given a topological space $X$, we denote by
$\C_c(X)$ the space of real-valued continuous functions on $X$ with
compact support.

Let $r>0$. A family $(\mu_{x})_{x\in \wt M}$ of nonzero finite
measures on $\partial_{\infty}\wt M$ whose support is the limit set
$\Lambda\Ga$ is a {\em Patterson density of dimension} $r$ for the
group $\Ga$ if it is $\Ga$-equivariant, that is, if it satisfies
\begin{equation}\label{eq:equivardensity}
\ga_*\mu_x=\mu_{\ga x}
\end{equation}
for all $\ga\in \Ga$ and $x\in\wt M$, and if the pairwise
Radon-Nikodym derivatives of the measures $\mu_x$ for $x\in \wt M$
exist and satisfy
\begin{equation}\label{eq:RNdensity}
\frac{d\mu_{x}}{d\mu_{y}}(\xi)=e^{-r\beta_{\xi}(x,y)}
\end{equation}
for all $x,y\in\wt M$ and  $\xi\in\partial_{\infty}\wt M$.

The {\em critical exponent} of $\Ga$ is
$$
\delta_{\Ga}=\lim_{n\to+\infty}\frac 1n\log\card\{\ga\in\Ga:
d(x_{0},\ga x_{0})\leq n\}. 
$$
The above limit exists and is positive, see \cite{Roblin02}, and the
critical exponent is independent of the base point $x_{0}$ used in its
definition. We assume that $\delta_{\Ga}$ is finite, which is in
particular the case if $M$ has a finite lower bound on its sectional
curvatures (see for instance \cite{Bowditch95}). We say that $\Ga$ is
{\it of divergence type} if its {\it Poincar\'e series}
$P_\Ga(s)=\sum_{\ga\in\Ga}e^{-sd(x_0,\ga x_0)}$ diverges at
$s=\delta_\Ga$.

Let $(\mu_{x})_{x\in \wt M}$ be a Patterson density of dimension
$\delta_{\Ga}$ for $\Ga$.  The {\em Bowen-Margulis measure} $\wt
m_{\rm BM}$ for $\Ga$ on $T^1\wt M$ is defined, using the Hopf
parametrisation, by
$$
d\wt m_{\rm BM}(v)=\frac{d\mu_{x_{0}}(v_{-})d\mu_{x_{0}}(v_{+})dt}
     {d_{x_{0}}(v_{-},v_{+})^{2\delta_{\Ga}}}
=e^{-\delta_{\Ga}(\beta_{v_{-}}(\pi(v),\,x_{0})+
\beta_{v_{+}}(\pi(v),\,x_{0}))}d\mu_{x_{0}}(v_{-})d\mu_{x_{0}}(v_{+})dt\,.    
$$
The Bowen-Margulis measure is independent of the base point $x_{0}$,
and its support is (in the Hopf parametrisation) $(\Lambda\Ga\times
\Lambda\Ga-\Delta)\times\RR$, where $\Delta$ is the diagonal in
$\Lambda\Ga\times\Lambda\Ga$.  It is invariant under the geodesic flow
and the action of $\Ga$, and thus it defines a measure $m_{\rm BM}$ on
$T^1M$, invariant under the quotient geodesic flow.  When the
Bowen-Margulis measure $m_{\rm BM}$ is finite, there exists a unique
(up to a multiplicative constant) Patterson density of dimension
$\delta_{\Ga}$, 
and the set of points in $T^1\wt M$ fixed by a nontrivial element of
$\Ga$ has measure $0$ for $\wt m_{\rm BM}$, see for instance
\cite[p.~19]{Roblin03}.  Denoting the total mass of a measure $m$ by
$\|m\|$, the probability measure $\frac{m_{\rm BM}}{\|m_{\rm BM}\|}$
is then uniquely defined. We will often make the assumption that
$m_{\rm BM}$ is finite, see the introduction for examples.

\medskip Let $C$ be a nonempty proper closed convex subset of $\wt
M$. We define the {\em skinning measure} $\wt\sigma_{C}$ of $\Ga$ on
$\normal C$, using the homeomorphism $w\mapsto w_+$ from $\normal C$
to $
\partial_{\infty}\wt M-\partial_{\infty}C$, by
\begin{align}
d\wt\sigma_{C}(w) &=  e^{-\delta_{\Ga}\beta_{w(+\infty)}(\pi(w),\,x_{0})}\,d(\nu
P_{C})_{*}(\mu_{x_{0}}|_{\partial_{\infty}\wt 
  M-\partial_{\infty}C})(w)\nonumber\\
  & =  e^{-\delta_{\Ga}\beta_{w_{+}}(P_C(w_+),\,x_{0})}\,d\mu_{x_{0}}(w_{+})\,.
\label{eq:defiskinmeas}
\end{align}
We will also consider $\wt\sigma_{C}$ as a measure on $T^1\wt M$ with
support contained in $\normal C$.

The skinning measure has been first defined by Oh and Shah \cite[\S
1.4]{OhShaCounting} for the outer unit normal bundles of spheres,
horospheres and totally geodesic subspaces in real hyperbolic spaces,
see also \cite[Lemma 4.3]{HerPau10} for a closely related measure. The
terminology comes from McMullen's proof of the contraction of the
skinning map (capturing boundary information for surface subgroups of
$3$-manifold groups) introduced by Thurston to prove his
hyperbolisation theorem.

When $C$ is a singleton $\{x\}$, we immediately have
$$
d\wt \sigma_C(w)=d\mu_x(w_+)\;.
$$
Let $w\in T^1\wt M$. When $C=H\!B_-(w)$ is the unstable horoball of
$w$, the measure 
$$
\musu{w}=\wt \sigma_{H\!B_-(w)}
$$
is the well known conditional measure of the Bowen-Margulis measure on
the strong unstable leaf $W^{\rm su}(w)$ of $w$ (see for instance
\cite{Margulis04,Roblin03}). Similarly, we denote by
$$
\muss{w}=\iota_*\big(\wt \sigma_{H\!B_{+}(w)}\big)
=\iota_*\big(\musu{\iota\,  w}\big)
$$
the conditional measure of the Bowen-Margulis measure on the strong
stable leaf $W^{\rm ss}(w)$ of $w$. These two measures are independent
of the element $w$ of a given strong unstable leaf and given strong
stable leaf respectively. For future use, using the homeomorphism
$v\mapsto v_-$ from $W^{\rm ss}(w)$ to $\partial_\infty\wt M-\{w_+\}$,
we have
\begin{equation}\label{eq:condtionelless}
d\muss{w}(v)=
e^{-\delta_{\Ga}\beta_{v_{-}}(P_{H\!B_+(w)}(v_-),\,x_{0})}\,d\mu_{x_{0}}(v_{-})\;.
\end{equation}
We also define the conditional measure of the Bowen-Margulis measure
on the stable leaf $W^{\rm s}(w)$ of $w$, using the homeomorphism
$(v',t)\mapsto v=\flow t v'$ from $W^{\rm ss} (w) \times\RR$ to
$W^{\rm s}(w)$, by
\begin{equation}\label{eq:defimus}
d\mus{w}(v)=e^{-\delta_\Ga t}\;d\muss{w}(v')dt\;.
\end{equation}
See for instance the assertion (iii) of the next proposition for an
explanation of the factor $e^{-\delta_\Ga t}$. We will not need in
this paper the similarly defined measure $d\muu{w}(v)= e^{\delta_\Ga
  t} \;d\musu{w}(v')dt$ on the unstable leaf $W^{\rm u}(w)$ of $w$.

\medskip
The following propositions collect some basic properties of the
skinning measures.

\bprop\label{prop:basics}
Let $C$ be a nonempty proper closed convex subset of $\wt M$, and let
$\wt\sigma_{C}$ be the skinning measure of $\Ga$ on $\normal C$. 

\noindent (i) The skinning measure $\wt\sigma_{C}$ is independent of
the base point $x_{0}$.

\noindent (ii) For all $\ga\in\Ga$, we have $\ga_*\wt\sigma_C=
\wt\sigma_{\ga C}$. In particular, the measure $\wt\sigma_{C}$ is
invariant under the stabiliser of $C$ in $\Ga$.

\noindent (iii) For all $s\geq 0$ and $w\in\normal C$, we have
$$
(\flow{s})_{*}\,\wt\sigma_{C}=e^{-\delta_\Ga s}\;\wt\sigma_{\N_{s}C} \;.
$$ 

\noindent (iv) The support of $\wt\sigma_{C}$ is $\{w\in\normal C:
w_{+}\in\Lambda\Gamma\}= \nu P_C(\Lambda\Gamma- \Lambda\Gamma
\cap \partial_\infty C)$.  In particular, $\wt \sigma_C$ is the zero
measure if and only if $\Lambda\Ga$ is contained in $\partial_\infty
C$. \eprop

It follows from (ii) that, for all $\ga\in\Ga$, we have
\begin{equation}\label{eq:invarcondstab}
\ga_*\musu{w}=\musu{\ga w},\;\;\;\ga_*\muss{w}=\muss{\ga w},
\;\;\;\ga_*\mus{w}=\mus{\ga w}\;.
\end{equation}
It follows from (iii) and from the equality $\iota\circ
\flow{t}=\flow{-t}\circ\iota$ that, for all $t\in\RR$, we have
\begin{equation}\label{eq:invargeoflocondstab}
(\flow t)_*\,\musu{w}=e^{-\delta_\Ga t}\;\musu{\flow t w},\;\;\;
(\flow {-t})_*\,\muss{w}=e^{-\delta_\Ga t}\;\muss{\flow{-t}  w},\;\;\;
(\flow t)_*\,\mus{w}=e^{\delta_\Ga t}\;\mus{w}\;.
\end{equation}

\dem The first assertion follows from Equation \eqref{eq:RNdensity}
with $r=\delta_\Ga$ and the second part of Equation
\eqref{eq:cocycle}.  The second assertion follows from Equation
\eqref{eq:equivardensity}, the first part of Equation
\eqref{eq:cocycle}, and the first assertion.

To prove the third assertion, we note that since $(\flow s w)_+=w_+$
and by the cocycle property \eqref{eq:cocycle}, we have 
\begin{align*}
  d\,\wt\sigma_{\N_{s}C}(\flow s w)&=
  e^{-\delta_{\Ga}\beta_{w_{+}}(\pi(\flow s
    w),\,x_{0})}\,d\mu_{x_{0}}(w_{+})=
  e^{-\delta_{\Ga}\beta_{w_{+}}(\pi(\flow s
    w),\,\pi(w))}\,d\,\wt\sigma_{C}( w)\\ &=e^{\delta_\Ga
    s}\,d\,\wt\sigma_{C}( w)\;.
\end{align*}

The fourth assertion follows from the fact that the support of any
Patterson measure is the limit set of $\Ga$. 
\cqfd

\medskip 
Given two nonempty closed convex subsets $C$ and $C'$ of $\wt
M$, let $\Omega_{C,C'}= \partial_{\infty}\wt M-(\partial_{\infty} C
\cup \partial_{\infty}C')$ and let
$$
h_{C,C'}:\nu P_{C}(\Omega_{C,C'})\to \nu P_{C'}(\Omega_{C,C'})
$$ 
be the restriction of $\nu P_{C'}\circ\nu P_{C}^{-1}$ to $\nu
P_{C}(\Omega_{C,C'})$. It is a homeomorphism between open subsets of
$\normal C$ and $\normal C'$, associating to the element $w$ in the
domain the unique element $w'$ in the range with $w'_+=w_+$.

\bprop\label{prop:abscontskinmeas} 
Let $C$ and $C'$ be nonempty proper closed convex subsets of $\wt M$
and let $h=h_{C,C'}$.  The measures $h_{*}\wt \sigma_{C}$ and $\wt
\sigma_{C'}$ on $\nu P_{C'}(\Omega_{C,C'})$ are absolutely continuous
one with respect to the other, with
$$
\frac{d\,h_{*}\wt\sigma_{C}}{d\,\wt\sigma_{C'}}(w')=
e^{\delta_{\Ga}\beta_{w_{+}}(\pi(w),\,\pi(w'))},
$$
for all $w\in \nu P_{C}(\Omega_{C,C'})$ and $w'=h(w)$.  
\eprop

\dem Since $w'_+=w_+$ and by the cocycle property \eqref{eq:cocycle},
we have
$$
d\wt\sigma_{C'}(w') =  
e^{-\delta_{\Ga}\beta_{w'_{+}}(P_{C'}(w'_+),\,x_{0})}\,
d\mu_{x_{0}}(w'_{+})
=e^{-\delta_{\Ga}\beta_{w_{+}}(P_{C'}(w'_+),P_{C}(w_+))}\,
d\wt\sigma_{C}(w)\;.
$$
Since $w=\nu P_C(w_+)$ and $\pi\circ\nu P_C= P_C$, and similarly for
$w'$, the result follows from the anti-symmetry of the Busemann
cocycle. \cqfd

\medskip We endow the set $\operatorname{Convex}(\wt M)$ of nonempty
closed convex subsets of $\wt M$ with the (metrisable, locally
compact) topology of the Hausdorff convergence on compact subsets: a
sequence $(C_{i})_{i\in\NN}$ of closed subsets of $\wt M$ converges to
a closed subset $C$ of $\wt M$ if and only if for every compact subset
$K$ in $\wt M$, the Hausdorff distance between $(C_{i}\cap K) \cup \,
^cK$ and $(C\cap K)\cup\;^cK$ tends to $0$. Note that being convex is
indeed a closed condition.  We endow the set $\operatorname{Measure}
(T^1 \wt M)$ of nonnegative regular Borel measures on $T^1\wt M$ with
the (metrisable, locally compact) topology of the weak-star
convergence: a sequence $(\mu_{i})_{i\in\NN}$ of such measures on $\wt
M$ converges to such a measure $\mu$ on $\wt M$ if and only if for
every compactly supported continuous function $f$ on $\wt M$, the
sequence $(\mu_{i}(f))_{i\in\NN}$ converges to $\mu(f)$.

\bprop\label{prop:skinningconv} The map from
$\operatorname{Convex}(\wt M)$ to $\operatorname{Measure}(T^1\wt M)$
which associates to $C$ its skinning measure $\wt \sigma_C$ is
continuous.  
\eprop

In particular, as the horoballs $H\!B_+(w)$ and $H\!B_-(w)$ depend
continuously on $w\in T^1\wt M$, the measures $\musu{w}$, $\muss{w}$
and $\mus{w}$ depend continuously on $w$.

\medskip
\dem Let $(C_{i})_{i\in\NN}$ be a sequence of nonempty closed convex
subsets of $\wt M$ which converges to a nonempty closed convex subset
$C$ for the Hausdorff convergence on compact subsets of $\wt M$, and
let us prove that $\wt\sigma_{C_{i}}\weakstar \wt \sigma_{C}$.

The sequence $(\normal C_{i})_{i\in\NN}$ of closed subsets of $T^1\wt
M$ converges to $\normal C$ for the Hausdorff convergence on compact
subsets of $T^1\wt M$. The sequence $(\partial_\infty\wt
M- \partial_\infty C_{i})_{i\in\NN}$ of open subsets of
$\partial_\infty\wt M$ converges to $\partial_\infty\wt
M- \partial_\infty C$ for the Caratheodory convergence (that is, for
the Hausdorff convergence of their complements).  Hence, the sequences
of maps $(P_{C_{i}})_{i\in\NN}$ and $(\nu P_{C_{i}})_{i\in\NN}$
converge to $P_{C}$ and $\nu P_{C}$ respectively for the uniform
convergence of maps on compact subsets $\partial_\infty\wt
M- \partial_\infty C$. Given two compact metric spaces $X$ and $Y$ and
a finite Borel measure $\mu$ on $X$, the pushforward map $f\mapsto
f_{*} \mu$ from the space of continuous maps from $X$ to $Y$ with the
uniform topology to the space of finite Borel measures on $Y$ with the
weak-star topology is continuous. The claim follows from these
observations, since the skinning measure on $C$ is a multiple by a map
depending continuously on $C$ of the pushforward by a map depending
continuously on $C$ of the fixed measure $\mu_{x_0}$.  \cqfd

\medskip The following result will be useful in Section
\ref{sec:equidistribution}.  Recall (see Equation
\eqref{eq:defiboulehamen}) that $V_{w,R}$ is the open ball of radius
$R$ and center $w$ in the strong stable leaf $W^{\rm ss}(w)$ of $w\in
T^1\wt M$ for the Hamenst\"adt distance.

\blemm\label{lem:veryoldlemseven} For every nonempty proper closed
convex subset $C$ in $\wt M$, there exists $R_0>0$ such that for every
$R\geq R_0$ and every $w\in\normal C$, we have
$\muss{w}(V_{w,R})>0$. If $\partial_{\infty} C\cap\Lambda\Ga \ne
\emptyset$, we may take $R_0=2$.  
\elemm

\medskip\noindent
\begin{minipage}{10.5cm} 
  \dem For all $w\in \normal C$ and $x\in C\cup \partial_{\infty} C$,
  by a standard comparison and convexity argument in the
  CAT$(-1)$-space $\wt M$ applied to the geodesic triangle with
  vertices $\pi(w), w_+,x$, the point $\pi(w)$ is at distance at most
  $2\log(\frac{1+\sqrt 5}{2})$ from the intersection between the
  stable horosphere $H_+(w)$ and the geodesic ray or line between $x$
  and $w_+$. Hence, by Equation \eqref{eq:majodistHamen}, for every
  $\xi'\in \partial_{\infty} C$, we have
$$
d_{W^{\rm ss}(w)}(w,\iota\,\nu
P_{W^{\rm ss}(w)}(\xi'))\leq \frac{1+\sqrt 5}{2}\;.
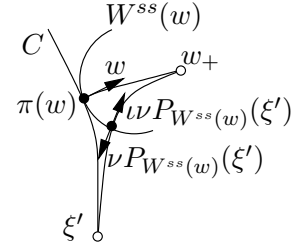
$$
\end{minipage}
\begin{minipage}{4.4cm}
\begin{center}
\input{fig_debexistR.pstex_t}
\end{center}
\end{minipage}

\medskip \noindent 
Thus, if $\partial_{\infty} C\cap\Lambda\Ga\ne\emptyset$, then we may
take $R_0 = 2 >\frac{1+\sqrt 5}{2}$, since by Proposition
\ref{prop:basics} (iv), the support of $\muss{w}$ is $\iota\,\nu
P_{HB_+(w)}(\Lambda\Gamma- \Lambda\Gamma\cap \{w^+\})$.

Assume now that $\partial_{\infty} C\cap\Lambda\Ga=\emptyset$. By
absurd, assume that, for all $n\in\NN$, there exists $w_n\in\normal C$
such that $\muss{w_n}(V_{w_n,n})=0$.  Assume first that
$(w_n)_{n\in\NN}$ has a convergent subsequence with limit $w\in
\normal C$. Since the measure $\muss{v}$ depends continuously on $v$,
for every compact subset $K$ of $W^{\rm ss}(w)$, we have
$\muss{w}(K)=0$. By Proposition \ref{prop:basics} (iv) and by Equation
\eqref{eq:condtionelless}, the support of the Patterson measure
$\mu_{x_0}$, which is the limit set of $\Ga$, is contained in
$\{w_{+}\}$. This is impossible, since $\Ga$ is nonelementary.

In the remaining case, the points $\pi(w_n)$ in $C$ converge, up to
extracting a subsequence, to a point $\xi$ in $\partial_{\infty}C$.
By definition of the map $\nu P_C$ and of $\normal C$, the points at
infinity $(w_n)_{+}$ converge to $\xi$. For every $\eta$ in
$\partial_\infty \wt M$ different from $\xi$, the geodesic lines from
$\eta$ to $(w_n)_{+}$ converge to the geodesic line from $\eta$ to
$\xi$. 

\medskip\noindent
\begin{minipage}{9cm} 
  ~~~ By convexity, if $n$ is big enough, the geodesic line
  $]\eta,(w_n)_{+}[$ meets $\N_1C$, hence passes at distance at most
  $2$ from $\pi(w_n)$. This implies (using Equation
  \eqref{eq:majodistHamen} as above) that if $n$ is big enough, then
  there exists $v\in V_{w_n, \,n}$ such that $\eta = v_-$.
\end{minipage}
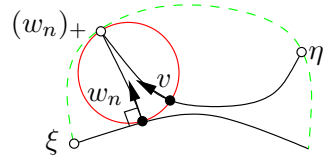
\begin{minipage}{5.9cm}
\begin{center}
\input{fig_existenceR.pstex_t}
\end{center}
\end{minipage}

\medskip \noindent Since we assumed that $\muss{w_n}(V_{w_n,n})=0$ for
all $n\in\NN$, Proposition \ref{prop:basics} (iv) implies that we
have $\eta\notin\Lambda\Ga$. Hence $\Lambda\Ga$ is contained in
$\{\xi\}$, a contradiction since $\Ga$ is nonelementary.  \cqfd

\medskip Let $C$ be a nonempty closed convex subset of $\wt M$, and
let $U_C$ be the open subset of $T^1\wt M$ defined in Equation
\eqref{eq:defiUC}. Note that $U_C$ has full Bowen-Margulis measure in
$T^1\wt M$ if the Patterson measure $\mu_{x}(\partial_\infty C)$ of
$\partial_\infty C$ is equal to $0$ (this being independent of $x\in\wt
M$), by the quasi-product structure of $\wt m_{\rm BM}$.

The following disintegration result of the Bowen-Margulis measure over
the skinning measure of $C$ is the crucial tool for the
equidistribution result in Section \ref{sec:equidistribution}. When
$\wt M$ has constant curvature and $\Ga$ is torsion free, this result
is implicit in \cite{OhShaCounting}.

\bprop\label{prop:disintegration} Let $C$ be a nonempty closed convex
subset of $\wt M$. The restriction to $U_C$ of the Bowen-Margulis
measure $\wt m_{\rm BM}$ disintegrates by the fibration $f_C:U_C\ra
\normal C$, over the skinning measure $\wt\sigma_{C}$ of $C$, with
conditional measure 
$e^{\delta_\Ga \beta_{w_+}(\pi(w),\,\pi(v))}\;d\mu^{\rm s} _w (v)$
on the fiber $f_C^{-1}(w)=W^s(w)$ of $w\in \normal
C$:
$$
d\wt m_{\rm BM}(v)=
\int_{w\in \normal C}e^{\delta_\Ga \beta_{w_+}(\pi(w),\,\pi(v))}\;
d\mu^{\rm s}_w(v)\;d\wt\sigma_{C}(w)\;.
$$
\eprop

\dem For every $\varphi\in\C_c(U_c)$, let $I_\varphi=\int_{v\in U_C}
\varphi(v)\;d\wt m_{\rm BM}(v)$. By the definition of $U_C$ and of the
Bowen-Margulis measure in the Hopf parametrisation, we have
$$
I_\varphi=\int_{v_+\in\,\partial_\infty\wt M-\partial_\infty C}
\int_{v_-\in\,\partial_\infty\wt M-\{v_+\}}\int_{t\in\RR}\varphi(v)\;
\frac{dt\,d\mu_{x_0}(v_{-})\, d\mu_{x_0}(v_{+})}
{d_{x_0}(v_-,v_+)^{2\delta_\Ga}}\;.
$$
\noindent
\begin{minipage}{5.9cm} 
  For every $v\in U_C$, let $w=f_C(v)=\nu P_C(v_+)$ and let $s\in\RR$
  be such that $v'=\flow{-s}v$ belongs to the strong stable leaf
  $W^{\rm ss}(w)$ of $w$. Note that, with $t$ the time parameter of
  $v$ in the Hopf parametrisation, the number $t-s$ depends only on
  $v_+$ and $v_-=v'_-$.  
\end{minipage}
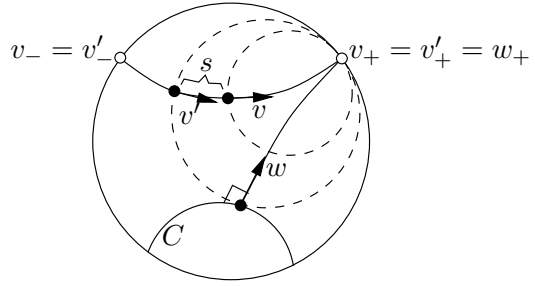
\begin{minipage}{9cm}
\begin{center}
\input{fig_desinteg.pstex_t}
\end{center}
\end{minipage}

\medskip \noindent Since the map from $\normal C$ to
  $\partial_\infty\wt M-\partial_\infty C$ defined by $w\mapsto w_+$
  and the map from $W^{\rm ss}(w)$ to $\partial_\infty\wt M-\{w_+\}$
  defined by $v'\mapsto v'_-$ are homeomorphisms, we have$$
I_\varphi=\int_{w\in\,\normal C}
\int_{v'\in\, W^{\rm ss}(w)}\int_{s\in\RR}\varphi(\flow s v')\;
\frac{ds\,d\mu_{x_0}(v'_{-})\, d\mu_{x_0}(w_{+})}
{d_{x_0}(v'_-,w_+)^{2\delta_\Ga}}\;.
$$
For every $w\in\normal C$ and $v'\in W^{\rm ss}(w)$, we claim
(explanations follow) that
\begin{align*}
\frac{d\mu_{x_0}(v'_{-})\, d\mu_{x_0}(w_{+})}
{d_{x_0}(v'_-,w_+)^{2\delta_\Ga}}&= 
\frac{e^{\delta_\Ga \beta_{v'_-}(\pi(v'),\,x_0))}\;
e^{\delta_\Ga \beta_{w_+}(\pi(w),\,x_0)}}
{e^{-\delta_\Ga(\beta_{v'_-}(x_0,\,\pi(v'))+
    \beta_{w_+}(x_0,\,\pi(v')))}}\; 
d\mu^{\rm  ss}_{w}(v')\,d\wt\sigma_C(w)\\&
= e^{\delta_\Ga\beta_{w_+}(\pi(w),\,\pi(v'))}\; 
d\mu^{\rm  ss}_{w}(v')\,d\wt\sigma_C(w)\\&
= 
d\mu^{\rm  ss}_{w}(v')\,d\wt\sigma_C(w)\;.
\end{align*}
The first equality holds by the definition of the measures $\mu^{\rm
  ss}_w$ (see Equation \eqref{eq:condtionelless}) and $\wt\sigma_{C}$
(see Equation \eqref{eq:defiskinmeas}), by the definition of the
visual distance $d_{x_0}$ (see Equation \eqref{eq:visdist}), and since
$\pi(v')$ belongs to the geodesic line between $v'_-$ and
$w_+=v'_+$. The second equality follows from the cocycle property
\eqref{eq:cocycle}. The third one holds since $\pi(w)$ and $\pi(v')$
both belong to the stable horosphere of $w$.

Hence, since $\beta_{w_+}(\pi(w), \,\pi(v))=s$ if $v=\flow s v'$ and
$v'\in W^{\rm ss}(w)$, and by the definition of the measure $\mu^{\rm
  s}_w$ (see Equation \eqref{eq:defimus}), we have
\begin{align}
I_\varphi&=\int_{w\in\,\normal C}
\int_{v'\in \,W^{\rm ss}(w)}\int_{s\in\RR}\varphi(\flow s v')\;
ds\,d\mu^{\rm  ss}_{w}(v')\,d\wt\sigma_C(w)
\label{eq:decomskinstrongstab}\\&
=\int_{w\in\,\normal C}\int_{v\in\, W^{\rm s}(w)}
\varphi(v)\;e^{\delta_\Ga\beta_{w_+}(\pi(w),\,\pi(v))}\;
d\mu^{\rm  s}_{w}(v)\,d\wt\sigma_C(w)\;,\nonumber
\end{align}
which proves the result.
\cqfd

\medskip We conclude this section by defining the skinning measures of
equivariant families of convex subsets.

Let $I$ be an index set endowed with a left action $(\ga,i)\mapsto \ga
i$ of $\Ga$.  A family $\D=(D_i)_{i\in I}$ of subsets of $\wt M$ or
$T^1\wt M$ indexed by $I$ is {\it $\Ga$-equivariant} if $\ga
D_i=D_{\ga i}$ for all $\ga\in\Ga$ and all $i\in I$.  We equip the
index set $I$ with the $\Ga$-equivariant equivalence relation $\sim$
(or $\sim_\D$ when we want to stress the dependence on $\D$), defined
by setting $i\sim j$ if and only if there exists $\ga\in\stab_\Ga D_i$
such that $j=\ga i$ (or equivalently if $D_j=D_i$ and $j=\ga i$ for
some $\ga\in\Ga$). Note that $\Ga$ acts on the left on the set of
equivalence classes $I/\!\!\sim$.

An example of such a family is given by fixing a subset $C$ of $\wt M$
or $T^1\wt M$, by setting $I=\Ga$ with the left action by translations
on the left $(\ga,i)\mapsto \ga i$, and by setting $D_i=i C$ for every
$i\in\Ga$.  In this case, we have $i\sim j$ if and only if $i^{-1}j$
belongs to the stabiliser $\Ga_C$ of $C$ in $\Ga$, and $I/\!\!\sim\;
=\Ga/\Ga_C$.  More general examples include $\Ga$-orbits of (usually
finite) collections of subsets of $\wt M$ or $T^1\wt M$ with (usually
finite) multiplicities.

A $\Ga$-equivariant family $(D_i)_{i\in I}$ of closed subsets of $\wt M$ or
$T^1\wt M$ is said to be {\it locally finite} if for every compact
subset $K$ in $\wt M$ or $T^1\wt M$, the quotient set $\{i\in I:
D_i\cap K\ne\emptyset\}/\!\sim $ is finite. In particular, the union
of the images of the sets $D_i$ by the map $\wt M\ra M$ or $T^1\wt
M\ra T^1M$ is closed. When $\Ga\bs I$ is finite, $(D_i)_{i\in I}$ is
locally finite if and only if, for all $i\in I$, the canonical map
from $\Ga_{D_i}\bs D_i$ to $M$ or $T^1M$ is proper, where $\Ga_{D_i}$
is the stabiliser of $D_i$ in $\Ga$.

Let $\D=(D_i)_{i\in I}$ be a locally finite $\Ga$-equivariant family
of nonempty proper closed convex subsets of $\wt M$. Then
$$
\wt\sigma_{\D}=\sum_{i\in I/\sim} \wt \sigma_{D_i}
$$
is a locally finite positive Borel measure on $T^1\wt M$ (independent
on the choice of representatives in $I/\!\!\sim$), called the {\em
  skinning measure} of $\D$ on $T^1\wt M$. It is $\Ga$-invariant by
Proposition \ref{prop:basics} (ii), and its support is contained in
$\bigcup_{i\in I/\sim} \normal D_i$. Hence $\wt\sigma_{\D}$ induces a
locally finite Borel positive measure $\sigma_{\D}$ on $T^1M=\Ga\bs
T^1\wt M$, called the {\em skinning measure} of $\D$ on $T^1M$.

For every $t\in[0,+\infty[\,$, let $\D_{t}=(\N_tD_i)_{i\in I}$, which
is also a $\Ga$-equivariant locally finite family of nonempty closed
convex subsets of $\wt M$. Note that by Proposition \ref{prop:basics}
(iii), we have
$$
(\flow t)_*\sigma_\D=e^{-\delta_\Ga t}\;\sigma_{\D_t}\;,
$$
and, in particular,
$$
\|\sigma_{\D_t}\|=e^{\delta_\Ga t}\;\|\sigma_{\D}\|\;.
$$
Note that the measure $\sigma_{\D_t}$ is finite if and only if the
measure $\sigma_{\D}$ is finite.

If the image in $M$ of the support of $\sigma_\D$ is compact, then
$\sigma_{\D}$ is finite. In particular, if $\Ga$ is geometrically
finite, the skinning measure of a Margulis neighbourhood of a cusp in
$\Ga\bs\wt M$ is finite, since for any parabolic fixed point $p$ of
$\Ga$, the quotient of $\Lambda\Ga-\{p\}$ by the stabiliser of $p$ in
$\Ga$ is compact.

Oh and Shah \cite[Theo.~1.5]{OhShaCounting} proved, in particular, that
$\|\sigma_\D\|$ is finite if $\Ga\bs I$ is finite, $\Ga$ is
torsion-free, $M$ is geometrically finite with constant curvature
$-1$, $\wt \D$ consists of codimension $1$ totally geodesic
submanifolds, and $\delta_\Ga>1$. See \cite[Theo.~6.4]{OhShaCounting}
for a statement without the codimension $1$ assumption, that we
generalise in the following section.

\medskip The next result relates the finiteness of the skinning measure
of $\D$ to the one of a nested family $\D'$.

\brema\label{rem:finitcontain} Let $\D=(D_i)_{i\in I}$ and
$\D'=(D'_i)_{i\in I}$ be locally finite $\Ga$-equivariant families of
nonempty proper closed convex subsets of $\wt M$, with $D'_i\subset
D_i$ for every $i\in I$.  Assume that

\noindent $\bullet$~ $P_{D'_i}(\xi)$ is the closest point in $D'_i$ to
$P_{D_i}(\xi)$, for every $\xi\in\Lambda \Ga-\partial_\infty D_i$;

\noindent $\bullet$~ there exists $c>0$ such that $d(P_{D_i}(\xi),
D'_i) \leq c$, for every $i\in I$ and $\xi\in \Lambda
\Ga- \partial_\infty D_i$;

\noindent $\bullet$~ for every $i\in I$, we have $\mu_{x_0}(\partial_\infty
D_i- \partial_\infty D'_i)=0$.

\noindent Then $\sigma_\D$ is finite if and only if $\sigma_{\D'}$ is
finite.  
\erema

It follows from this remark that for every $\epsilon\geq 0$, if
$\D''=(\N_\epsilon D_i)_{i\in I}$, then $\sigma_{\D''}$ is finite if
and only if  $\sigma_{\D}$ is finite.

The first assumption is also satisfied if $\wt M$ has constant
curvature $-1$ and $D_i$ is totally geodesic for all $i\in I$, since
by homogeneity, for every $\xi$ in $\wt M$ and $x\neq y$ in $\wt M$
such that $\angle_x(\xi,y)= \frac{\pi}{2}$, the value $\beta_\xi(y,x)$
is a strictly increasing function of only $d(x,y)$.

\medskip \dem By the first assumption, the map $\theta:\nu
P_{D_i}(\Lambda\Ga-\partial_\infty D_i) \ra \nu
P_{D'_i}(\Lambda\Ga-\partial_\infty D'_i)$ defined by $w\mapsto w'$
where $w'_+=w_+$ and $\pi(w')$ is the closest point on $D'_i$ to
$\pi(w)$ is a homeomorphism onto is image such that
$$
{\nu P_{D'_i}}_{\mid \Lambda\Ga-\partial_\infty D_i}=
\theta \circ {\nu P_{D_i}}_{\mid \Lambda\Ga-\partial_\infty D_i}\;.
$$
By the definition of the skinning measures, using this homeomorphism
$\theta$, we have, for all $w'\in \theta(\nu
P_{D_i}(\Lambda\Ga-\partial_\infty D_i))$,
$$
d\,\wt\sigma_{D'_i}(w')=e^{-\delta_{\Ga}\beta_{w'_{+}}(P_{D'_i}(w'_+),\,P_{D_i}(w'_+))}
\,d\,\theta_*\wt\sigma_{D_i}(w')\;.
$$
The result then follows by the second and third assumptions.
\cqfd

\section{Finiteness and fluctuation of  the skinning measure}
\label{sect:finitskinn}

We will say that a discrete group $\Ga'$ of isometries of
$\wt M$ has {\it regular growth} if there exists $c>0$ such that for
every $N\in\NN$, we have
$$
\frac{1}{c}\;e^{\delta_{\Ga'}\,N}\leq \card\{\ga\in\Ga'\;:\;d(x_0,\ga x_0) 
\leq N\}\leq c\;e^{\delta_{\Ga'}\,N}\;.
$$
This does not depend on $x_0$, and the upper bound holds for all
nonelementary groups $\Ga'$ (see for instance \cite[page
11]{Roblin03}). If the Bowen-Margulis measure $m_{\rm BM}$ on $T^1M$
is finite, then $\Ga$ has regular growth (there even exists an
explicit $c>0$ such that $\card\{\ga\in\Ga\;:\;d(x_0,\ga x_0) \leq
N\}\sim c\;e^{\delta_{\Ga}N}$, see for instance \cite{Roblin03}).  If
$\wt M$ is a symmetric space, then any discrete parabolic group of
isometries of $\wt M$ has regular growth. In particular, if $\wt M$ is
the real hyperbolic space $\HH^n_\RR$, then by a theorem of
Bieberbach, any discrete parabolic group $\Ga'$ contains a finite
index subgroup isomorphic to $\ZZ^k$ for some $k\in \{0,\dots, n-1\}$
called the {\it rank} of the fixed point of $\Ga'$, and an easy and
well known computation in hyperbolic geometry proves that the critical
exponent of $\Ga'$ is
\begin{equation}\label{eq:expocritparahypreal}
\delta_{\Ga'}=\frac{k}{2}\;,
\end{equation} 
and that $\Ga'$ has regular growth. Note that there exist complete
simply connected Riemannian manifolds with pinched negative curvature
having discrete parabolic groups of isometries which do not have
regular growth, see for instance \cite{DalOtaPei00}.

We will say that a convex subset $C$ of $\wt M$ is {\it almost
  cone-like in cusps} for a discrete group $\Ga'$ of isometries of
$\wt M$ if for any parabolic point $p'$ of $\Ga'$ belonging to
$\partial_\infty C$ and any horoball $\H'$ centered at $p'$, there
exists $r\geq 0$ and $x'_0\in\partial \H'$ such that $C\cap \H'\cap
\N_{2\log(1+\sqrt{2})}(\C\Lambda\Ga')$ is contained in the orbit of
$\N_r([x'_0,p'[)$ under the stabiliser in $\Ga'$ of $p$ and $C$. It
follows from the arguments of \cite[\S 4]{OhShaCounting} that if $\wt
M$ has constant sectional curvature $-1$, if $C$ is a totally geodesic
submanifold and if $\Ga'$ is torsion free and geometrically finite,
then $C$ is almost cone-like in cusps for $\Ga'$.

\btheo\label{theo:geneohshah} Let $\wt M$ be a complete simply
connected Riemannian manifold with sectional curvature bounded above
by $-1$. Let $\Ga$ be a geometrically finite discrete group of
isometries of $\wt M$, of divergence type, with finite critical
exponent.  Let $\D=(D_i)_{i\in I}$ be a locally finite
$\Ga$-equivariant family of nonempty proper convex subsets of $\wt M$
which are almost cone-like in cusps for $\Ga$, with $\Ga\bs I$ finite.
Assume that for every parabolic point $p$ of $\Ga$ and every $i\in I$
such that $p\in\partial_\infty D_i$, the stabilisers $\Ga_p$ and
$\Ga_{D_i}$ in $\Ga$ of $p$ and $D_i$, respectively, have regular
growth and satisfy
\begin{equation}\label{eq:condiexpocrit}
\delta_\Ga>\;2(\delta_{\Ga_p}-\delta_{\Ga_{D_i}\cap\Ga_p})\;.
\end{equation} 
Then the skinning measure $\sigma_{\D}$ of $\D$ on $T^1M$ is finite.
\etheo

We make some comments on this statement before giving its proof.

\medskip
\noindent{\bf Remarks.} (1) When $\wt M$ is a symmetric space (in
particular when $\wt M$ has constant sec\-tional curvature $-1$),
every geometrically finite group of isometries of $\wt M$ is of
divergence type. This is not true in general, but holds true if
$\delta_\Ga> \delta_{\Ga_p}$ for every parabolic point $p$ of $\Ga$,
see \cite{DalOtaPei00}.  As already said, $\delta_{\Ga}$ is finite for
instance if $M$ has a finite lower bound on its sectional curvatures.

(2) Assume in this remark that the index of $\Ga_{D_i}\cap\Ga_p$ in
$\Ga_p$ is finite for every parabolic point $p$ of $\Ga$ and every
$i\in I$ such that $p\in\partial_\infty D_i$. Then $\delta_{\Ga_p}=
\delta_{\Ga_{D_i}\cap\Ga_p}$, and this equality implies that the
condition \eqref{eq:condiexpocrit} is satisfied. When $\wt M$ has
constant sectional curvature $-1$, the subsets $D_i$ are totally
geodesic submanifolds, and $\Ga$ is torsion-free, the finiteness of
$\sigma_{\D}$ follows from \cite[Theo.~6.3]{OhShaCounting}.

(3) Assume in this remark that $\wt M$ has constant sectional
curvature $-1$ and that the subsets $D_i$ are totally geodesic
submanifolds. Let us prove that for every parabolic point $p$ of $\Ga$
belonging to $\partial_\infty D_i$, we have $\delta_{\Ga_p}-
\delta_{\Ga_{D_i}\cap\Ga_p} \leq \frac{1}{2}\operatorname{codim}(D_i)$
(see also \cite[Lem.~6.2]{OhShaCounting} when $\Ga$ is torsion-free).
This will imply that the condition \eqref{eq:condiexpocrit} is
satisfied if $\delta_\Ga>1$ and if the elements of $\D$ have
codimension $1$.

Let $k$ be the rank of $\Ga_p$. In particular, $\delta_{\Ga_p} =
\frac{k}{2}$ by Equation \eqref{eq:expocritparahypreal}.  Up to taking
a finite index subgroup, and choosing appropriate coordinates, we may
assume that $p$ is the point at infinity in the upper halfspace model
of $\wt M=\HH^n_\RR$, that $\Ga_p$ is the lattice $\ZZ^k$ of $\RR^k$
acting by translations on the first factor (and trivially on the
second one) on $\RR^k\times\RR^{n-k-1}=\RR^{n-1}
= \partial_\infty\HH^n_\RR-\{p\}$, and that $E=\partial_\infty
D_i-\{p\}$ is a linear subspace of $\RR^{n-1}$. Let $F=E\cap \RR^k$,
which is a linear subspace of $\RR^k$. Since the family $\D$ is
locally finite, the image of $F$ in the torus $\RR^k/\ZZ^k$ is closed,
hence it is  a subtorus. Since $0\in F$, the subgroup $\ZZ^k\cap F$ is
hence a lattice in $F$.  Therefore, by Equation
\eqref{eq:expocritparahypreal}, 
$$
2(\delta_{\Ga_p}-
\delta_{\Ga_{D_i}\cap\Ga_p})=\operatorname{codim}_{\RR^k}(F)\leq
\operatorname{codim}_{\RR^{n-1}}(E)=\operatorname{codim}(D_i)\;.
$$

(4) Theorem \ref{theo:geneohshah} is optimal, since when $\wt M$ has
constant sec\-tional curvature $-1$, the subsets $D_i$ are totally
geodesic submanifolds and $\Ga$ is torsion-free, it is proved in
\cite[Theo.~6.4]{OhShaCounting} that the validity of Equation
\eqref{eq:condiexpocrit} (translated using Equation
\eqref{eq:expocritparahypreal}), for all $i,p$ as in the statement, is
a necessary and sufficient condition for the skinning measure
$\sigma_{\D}$ to be finite.

(5) The ideas of the proof of Theorem \ref{theo:geneohshah} are a
blend of the ones of the finiteness of the Bowen-Margulis measure
under a separation condition on the critical exponents in
\cite{DalOtaPei00} (see also \cite{PauPolSha11} for the case of Gibbs
measures), and the ones of a generalisation to variable curvature of
Sullivan's fluctuating density property in \cite[\S 4]{HerPau04}.

\medskip \dem We may assume that $\Ga\bs I$ is a singleton. Let us fix
$i\in I$. We may assume that $\partial_\infty D_i\cap \Lambda\Ga$ is
nonempty. Otherwise indeed, since $\nu P_{C_i}$ is a homeomorphism and
$\Lambda\Ga$ is closed, the support of $\wt\sigma_{D_i}$, which is the
set of elements $v\in\partial^1_+D_i$ such that $v_+\in \Lambda\Ga$
(see Proposition \ref{prop:basics} (iv)), is compact. Hence the
support of $\sigma_{\D}$ is compact, therefore $\sigma_{\D}$ is
finite. Let $\pi:T^1\wt M\ra \wt M$ and again $\pi:T^1 M\ra M$ be the
base point projections. Note that the skinning measure $\sigma_{\D}$ is
finite if and only if its pushforward measure $\pi_* \sigma_{\D}$ is
finite.

In what follows, let $\varepsilon=\ln(1+\sqrt{2})$: Note that
for any geodesic triangle in $\HH^2_\RR$ with two ideal vertices and a
right angle at the  vertex $x\in\HH^2_\RR$, the distance from $x$ to its
opposite side is exactly $\varepsilon$.

\blemm \label{lem:supportalmostconvhull}
The support of the measure $\pi_*\wt \sigma_{D_i}$, which is
$\{P_{D_i}(\xi)\;:\; \xi\in\Lambda\Ga-\partial_\infty D_i\}$, is
contained in the closed $\varepsilon$-neighbourhood of the convex hull
$\C\Lambda\Ga$.  
\elemm

\dem Let $\xi\in\Lambda\Ga-\partial_\infty D_i$, let
$\xi'\in \partial_\infty D_i\cap \Lambda\Ga$, and let $x$ be the
closest point to $\xi$ on $D_i$. Then the geodesic ray from $x$ to
$\xi'$, which is contained in $D_i$ by convexity, makes an angle at
least $\frac{\pi}{2}$ at $x$ with the geodesic ray from $x$ to
$\xi$. By a standard comparison result and the definition of
$\varepsilon$, the point $x$ is hence at distance at most
$\varepsilon$ from the geodesic line between $\xi$ and $\xi'$, which
is contained in $\C\Lambda\Ga$.  \cqfd

\medskip Let $\Par_\Ga$ be the set of parabolic fixed points of
$\Ga$. Since $\Ga$ is geometrically finite (see for instance
\cite{Bowditch95}),

\noindent $\bullet$~ every $p\in\Par_\Ga$ is {\it bounded}, that is,
its stabiliser $\Ga_p$ in $\Ga$ acts properly with compact quotient on
$\Lambda\Ga-\{p\}$;

\noindent $\bullet$~ the action of $\Ga$ on $\Par_\Ga$ has only
finitely many orbits;

\noindent $\bullet$~ there exists a $\Ga$-invariant family
$(\H_p)_{p\in\Par_\Ga}$ of pairwise disjoint closed horoballs, with
$\H_p$ centered at $p$, such that the quotient
$$
M_0=\Ga\bs\big(\C\Lambda\Ga-\bigcup_{p\in\Par_\Ga} \H_p\big)
$$ 
is compact. The inclusion $\H_p\subset \wt M$ induces an injection
$\Ga_p\bs\H_p\ra\Ga\bs \wt M$ and we will identify $\Ga_p\bs\H_p$ with
its image in $\Ga\bs\wt M$. In particular, $\H_p$ is {\it precisely
  invariant} under $\Ga$, that is, for all $\ga\in\Ga-\Ga_p$, we have
$\ga\H_p\cap \H_p=\emptyset$.

By Lemma \ref{lem:supportalmostconvhull} (and since the
$\varepsilon$-neighbourhood of $M_0$ is also compact), we hence only
have to prove the finiteness of $\pi_*\sigma_{\D}(\Ga_p\bs\H_p)$ for
all $p\in\Par_\Ga$. By the local finiteness of $\D$ and the fact that
parabolic fixed points are bounded, for all $p\in\Par_\Ga$, if the
orbit $\Ga p$ does not meet $\partial_\infty D_i$, then
$\pi_*\sigma_{\D}(\Ga_p\bs\H_p)$ is finite.

We hence assume that there exist $p\in\Par_\Ga
\cap p\in \partial_\infty D_i$, and we want to prove the finiteness of
$\pi_*\sigma_{\D}(\Ga_p\bs\H_p)$. To simplify the notation, let
$\Ga_{p,i}=\Ga_{D_i}\cap\Ga_p$, $\delta_{p,i}=\delta_{\Ga_{p,i}}$,
$\delta_{p}=\delta_{\Ga_{p}}$ and $\delta=\delta_\Ga$. Let $x_0$ be a
point in $D_i\cap \partial \H_p$ (which exists up to shrinking
$\H_p$). Since $p$ is the endpoint of a geodesic ray contained in
$D_i$ and of a geodesic ray contained in $\C\Lambda\Ga$, and since
geodesic rays with the same point at infinity become arbitrarily
close, up to shrinking $\H_p$, we may assume that $x_0\in
\N_\varepsilon (\C\Lambda\Ga)$.

Choose a set of representatives $\Ga_p\bs\!\bs \Ga$ of the right
cosets in $\Ga_p\bs \Ga$ such that for all $\ga'\in \Ga_p\bs\!\bs
\Ga$, we have
$$
d(x_0, \ga'x_0) = \min_{\alpha\in\Ga_p} \;d(x_0, \alpha\ga'x_0)\;.
$$

Choose a set of representatives $\Ga_{p,i}\bs\!\bs \Ga_p$ of the right
cosets in $\Ga_{p,i}\bs \Ga_p$ such that for all $\overline{\alpha}\in
\Ga_{p,i}\bs\!\bs \Ga_p$, we have
$$
d(x_0, \overline{\alpha}\,x_0) = 
\min_{\beta\in\Ga_{p,i}} \;d(x_0, \beta\overline{\alpha}\,x_0)\;.
$$
Note that any $\ga\in\Ga$ may be uniquely written $\ga=\beta\ov
\alpha\ga'$ with $\beta\in\Ga_{p,i}$, $\ov\alpha\in \Ga_{p,i}
\bs\!\bs\Ga_p$ and $\ga'\in\Ga_{p}\bs\!\bs\Ga$.

\blemm \label{lem:quasigeod} 
There exists $c_1>0$ such that the following assertions hold.

\noindent (i) For all $\ga'\in\Ga_p\bs\!\bs \Ga$, the closest point on
$\H_p$ to $\ga'x_0$ is at distance at most $c_1$ from
$x_0$. Furthermore, for all $\ga'\in\Ga_p\bs\!\bs \Ga$ and
$\alpha\in\Ga_p$, for every $y$ in the geodesic ray $[x_0,p[\,$, we
have
$$
d(y, \alpha y)+d(y, x_0)+d(x_0, \ga'x_0)-c_1\leq 
d(y, \alpha\ga'x_0)\leq d(y, \alpha y)+d(y, x_0)+d(x_0, \ga'x_0)\;.
$$

\noindent (ii) For all $\overline{\alpha}\in\Ga_{p,i}\bs\!\bs \Ga_p$,
the closest point on $D_i$ to $\ov \alpha x_0$ is at distance at most
$c_1$ from the geodesic ray $[x_0,p[$.  Furthermore, for all
$\overline{\alpha}\in\Ga_{p,i}\bs\!\bs \Ga_p$ and $\beta\in\Ga_{p,i}$,
$$
\max\{d(x_0, \overline{\alpha}x_0),\;d(x_0, \beta x_0)\}-c_1
\leq d(x_0, \beta\overline{\alpha}\,x_0)\leq 
\max\{d(x_0, \overline{\alpha}x_0),\;d(x_0, \beta x_0)\}+c_1\;.
$$
\elemm

\dem (i) For all $\ga\in\Ga$, let $p_\ga$ be the closest point to $\ga
x_0$ on $\H_p$, which lies on the geodesic ray $[\ga x_0,p[\,$. Hence,
by our choice of $x_0$ and by convexity, $p_\ga$ is at bounded
distance from $\C\Lambda\Ga$. Since $\H_p$ is precisely invariant and
$x_0\in\partial \H_p$, the point $p_\ga$ belongs to
$\partial\H_p$. For all $\ga\in\Ga$ and $\alpha\in\Ga_p$, if
$p_\ga\neq \alpha^{-1}x_0,\ga x_0$, then the angle at $p_\ga$ between
$[p_\ga, \alpha^{-1}x_0]$ and $[p_\ga, \ga x_0]$ is at least
$\frac{\pi}{2}$ by the convexity of $\H_p$. Hence, by a standard
comparison argument and the definition of $\varepsilon$, the distance
between $p_\ga$ and $[\alpha^{-1}x_0,\ga x_0]$ is at most
$\varepsilon$. By the triangle inequality, we have
$$
d(\alpha^{-1} x_0,p_\ga)+d(p_\ga, \ga x_0)-2\varepsilon\leq 
d(x_0, \alpha\ga x_0)\leq d(\alpha ^{-1} x_0,p_\ga)+d(p_\ga, \ga x_0)\;.
$$
These inequalities are also true if $p_\ga$ is equal to $\alpha^{-1}
x_0$ or to $\ga x_0$. Since $p_\ga$ is at bounded distance from
$\C\Lambda\Ga \cap \partial\H_p$ and since the action of $\Ga_p$ on
$\C\Lambda\Ga\cap\partial\H_p$ is cocompact, there exists
$\alpha_\ga\in\Ga_p$ such that $d(p_\ga, \alpha_\ga x_0)$ is bounded,
say by $c'_1$.  Let $\ga'\in \Ga_p\bs\!\bs \Ga$.  Assume for a
contradiction that $d(x_0,p_{\ga'})>2\varepsilon +c'_1$. Then, using
the above centered equation with $\alpha=1$ and $\ga=\ga'$, we have
\begin{align*}
d(\alpha_{\ga' }^{-1}\ga' x_0, x_0) & = d(\ga' x_0,\alpha_{\ga' }x_0)\le 
d(\ga' x_0, p_{\ga'}) +d(p_{\ga'},\alpha_{\ga' }x_0)\\ &
\leq d(\ga' x_0, x_0)- d(x_0, p_{\ga'}) +2\varepsilon +c'_1
 <d(\ga' x_0, x_0)\;,
\end{align*}
which contradicts the minimality property of $d(\ga' x_0, x_0)$. 
This proves the first claim of Assertion (i) if  $c_1\ge 2\epsilon + c'_1$.

The first claim and the convexity of the horoball of center $p$ whose
boundary contains $y$ (which implies that if $y\neq x_0,\alpha^{-1}y$,
then the angle at $y$ between $[y,x_0]$ and $[y,\alpha^{-1}y]$ is at
least $\frac{\pi}{2}$) imply that the length of the piecewise geodesic
$[\ga'x_0,x_0]\cup[x_0,y] \cup[y,\alpha^{-1}y]$ is almost additive,
yielding the left hand side of the second claim of Assertion (i).  Its
right hand side follows by the triangle inequality.

\medskip (ii) For all $\alpha\in\Ga_p$, let $q_\alpha$ be the closest
point to $\alpha x_0$ on $D_i$.  By the convexity of $\H_p$ and since
$\alpha x_0\in \partial \H_p$, we have $q_\alpha\in\H_p$. By the
convexity of $D_i$ and as in (i), the point $q_\alpha$ lies at
distance at most $\varepsilon$ of the geodesic ray $[\alpha
x_0,p[\,$. Since $x_0\in\N_\varepsilon (\C\Lambda\Ga)$, the point
$q_\alpha$ is at distance at most $2\varepsilon$ from a point in
$\C\Lambda\Ga$.  Hence, $q_\alpha\in D_i\cap \H_p\cap
\N_{2\varepsilon}(\C\Lambda\Ga)$. Since $D_i$ is almost cone-like
in cusps for $\Ga$, there exists $\beta_\alpha\in \Ga_{p,i}$ such that
the distance between $\beta_\alpha q_\alpha= q_{\beta_\alpha\alpha}$
and $[x_0,p[$ is less than a constant.

\medskip \noindent
\begin{minipage}{9.9cm} ~~~ Let $q'_\alpha$ be the closest point to
  $q_\alpha$ on $[x_0,p[$. By quasi-geodesic arguments, there exists a
  constant $c>0$ such that
$$
|\,d(x_0, \beta_\alpha \alpha x_0)- 2d(x_0,\beta_\alpha q_\alpha)\,|
\leq c\;,
$$
$$
|\,d(x_0, \alpha x_0)- 2d(x_0, \beta_\alpha q_\alpha)-2d(q_\alpha,
q'_\alpha)\,|\leq c\;.
$$
Using a similar argument to that used in the proof of Assertion (i),
this proves that $q_{\ov\alpha}$ is at distance less than a constant
from $[x_0,p[$ for every $\ov\alpha\in\Ga_{p,i}\bs\!\bs \Ga_p$.
\end{minipage}
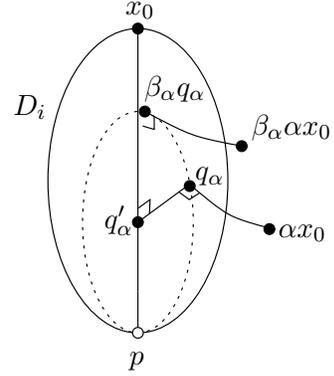
\begin{minipage}{5cm}
\begin{center}
\input{fig_projDidansHp.pstex_t}
\end{center}
\end{minipage}

\medskip\noindent 
For all $\ov\alpha\in\Ga_{p,i}\bs\!\bs \Ga_p$ and $\beta\in\Ga_{p,i}$,
since $\beta^{-1}x_0\in D_i$ and $q_{\ov\alpha}$ is the closest point
to $\ov\alpha x_0$ on $D_i$, we have
$$
d(\beta^{-1} x_0,q_{\ov\alpha})+d(q_{\ov\alpha}, \ov\alpha x_0)
-2\varepsilon\leq 
d(\beta^{-1} x_0,\ov\alpha x_0)\leq 
d(\beta^{-1} x_0,q_{\ov\alpha})+d(q_{\ov\alpha}, \ov\alpha x_0)\;.
$$

\noindent
\begin{minipage}{10.7cm} ~~~ For every $\alpha'\in\Ga_p$, let
  $r_{\alpha'}$ be the closest point to $\alpha'x_0$ on
  $[x_0,p[\,$. Hence by the above argument, there exists $c'>0$ such
  that
$$
|\,d(x_0, \beta\ov\alpha x_0)-d(\beta^{-1} x_0,r_{\ov\alpha})
-d(r_{\ov\alpha}, \ov\alpha x_0)\,|\leq c'\;.
$$
For all $y\in [x_0,p[$, we have 
$$
d(\alpha' x_0,r_{\alpha'})+
d(r_{\alpha'},y)-2\varepsilon\leq d(\alpha'
x_0,y)\leq d(\alpha' x_0,r_{\alpha'})+
d(r_{\alpha'},y)\;.
$$
\end{minipage}
\begin{minipage}{4.2cm}
\begin{center}
\input{fig_horoarbre.pstex_t}
\end{center}
\end{minipage}

\medskip\noindent
\begin{minipage}{10cm} ~~~ Let $H'$ be the horoball centered at $p$
  whose boundary contains $r_{\alpha'}$ and let $s$ be the
  intersection point of $[\alpha'x_0,p[$ with $\partial H'$.  Then
  $$
  d(\alpha' x_0,r_{\alpha'})\geq d(\alpha' x_0,s)=d(r_{\alpha'},x_0)\;,
  $$
  since $x_0$ and $\alpha' x_0$ are on the same horosphere centered at
  $p$. 
\end{minipage}
\begin{minipage}{4.9cm}
\begin{center}
\input{fig_projhoro.pstex_t}
\end{center}
\end{minipage}

\medskip By an easy comparison argument in the geodesic triangle with
vertices $r_{\alpha'}$, $\alpha'x_0$ and $p$, we have
$d(s,r_{\alpha'})\leq 1$. Hence
$$
d(\alpha' x_0,r_{\alpha'})\leq 
d(\alpha' x_0,s)+d(s,r_{\alpha'})\leq d(x_0,r_{\alpha'})+ 1\;.
$$

Applying this for $\alpha'=\beta^{-1},\ov\alpha$ and
$y=r_{\ov\alpha},r_{\beta^{-1}},x_0$, we have
$$
|\,d(x_0, \beta\ov\alpha x_0)-d(\beta^{-1} x_0,x_0)\,|\leq
c'+1+2\varepsilon
$$
if $r_{\ov\alpha}$ belongs to $[x_0,r_{\beta^{-1}}]$, and otherwise
$$
|\,d(x_0, \beta\ov\alpha x_0)-d(x_0,\ov\alpha x_0)\,|\leq
c'+1+2\varepsilon\;.
$$
This proves the result.
\cqfd

\medskip The next lemma, which uses the regular growth property of
$\Ga_p$ and $\Ga_{p,i}$, implies, in particular,  that the ``relative''
critical exponent of $\Ga_p$ modulo $\Ga_{p,i}$ is
$\delta_{p}-\delta_{p,i}$ (see for instance \cite{Paulin13b} for
background on relative Poincar\'e series).

\blemm\label{lem:expcritrelatif}
There exists $c_2>0$ such that for every $t\in[0,+\infty[\,$, we have
$$
\frac{1}{c_2}\;e^{(\delta_{p}-\delta_{p,i})t}\leq
\card\{\ov\alpha\in\Ga_{p,i}\bs\!\bs\Ga_p\;:\;d(x_0,\ov\alpha\, x_0)
\leq t\} \leq c_2\;e^{(\delta_{p}-\delta_{p,i})t}\;.
$$
\elemm

\dem For all $t\in[0,+\infty[\,$, define
$$
f(t)=\card\{\alpha\in\Ga_p\;:\;d(x_0,\alpha x_0)\leq t\}
\;\;\;{\rm and}\;\;\;
g(t)=\card\{\beta\in\Ga_{p,i}\;:\;d(x_0,\beta x_0)\leq t\}\;.
$$
Since $\Ga_p$ and $\Ga_{p,i}$ have regular growth, there exists a
constant $c>0$ such that for all $t\in[0,+\infty[$, we have
$$
\frac{1}{c}\,e^{\delta_p\, t}\leq f(t)\leq c\,e^{\delta_p\, t}
\;\;\;{\rm and}\;\;\;
\frac{1}{c}\,e^{\delta_{p,i} \,t}\leq g(t)\leq 
c\,e^{\delta_{p,i}\, t}\;.
$$
Also define $E=\Ga_p\times (\Ga_{p,i}\bs\!\bs\Ga_p)$ and $h(t)=
\card\{\ov\alpha\in\Ga_{p,i}\bs\!\bs\Ga_p\;:\; d(x_0,\ov\alpha x_0)\leq
t\}$.

For all $t\geq c_1$, we have, using Lemma \ref{lem:quasigeod} (ii) to
get the inequality,
\begin{align*}
  f(t-c_1)&=\card\{(\beta,\ov\alpha)\in E \;:\;d(x_0,\beta\ov\alpha
  x_0)\leq t-c_1\}\\& =\card\{(\beta,\ov\alpha)\in E
  \;:\;d(x_0,\beta\ov\alpha x_0)\leq t-c_1, \;d(x_0,\beta x_0)\leq
  d(x_0,\ov\alpha x_0) \}\\&\;\;\;\;+ \card\{(\beta,\ov\alpha)\in E
  \;:\;d(x_0,\beta\ov\alpha x_0)\leq t-c_1, \;d(x_0,\beta x_0)>
  d(x_0,\ov\alpha x_0) \}\\& \leq \card\{(\beta,\ov\alpha)\in E
  \;:\;d(x_0,\ov\alpha x_0)\leq t,\; d(x_0,\beta x_0)\leq t\}\\
  &\;\;\;\;+ \card\{(\beta,\ov\alpha)\in E \;:\;d(x_0,\beta
  x_0)\leq t, \;t\geq d(x_0,\ov\alpha x_0) \}\\
  & =2\;g(t)\,h(t)\;.
\end{align*}
This gives the lower bound in Lemma \ref{lem:expcritrelatif}.

Similarly, for all $t\geq c_1$, we have
\begin{align*}
  f(t+c_1+1)&\geq \card\{(\beta,\ov\alpha)\in E
  \;:\;\\&\;\;\;\;\;\;\;\; t-c_1<d(x_0,\beta\ov\alpha x_0)\leq
  t+c_1+1, \;d(x_0,\beta x_0)\leq d(x_0,\ov\alpha x_0) \}\\& \geq
  \card\{(\beta,\ov\alpha)\in E \;:\;t<d(x_0,\ov\alpha x_0)\leq t+1,
  \;d(x_0,\beta x_0)\leq t+1 \}\\
  & =g(t+1)(h(t+1)-h(t))\;.
\end{align*}
A geometric series summation argument gives the upper bound in Lemma
\ref{lem:expcritrelatif}.  
\cqfd

\medskip Now, let $\F^+_{p,i}$ be the set of accumulation points in
$\partial_\infty \wt M$ of the orbit points $\overline{\alpha}\ga'x_0$
where $\overline{\alpha} \in\Ga_{p,i}\bs\!\bs \Ga_p$ and $\ga'\in
\Ga_p\bs\!\bs \Ga$.  

\blemm \label{lem:eq:decomplimset} We have
$
\Lambda\Ga= \{p\}\cup
\bigcup_{\beta\in\Ga_{p,i}}\beta \F^+_{p,i}\;.
$
\elemm

Note that in general this union is not a disjoint union.  

\medskip \dem Every element $\xi$ in $\Lambda\Ga$ is the limit of a
sequence $(\beta_i \overline{\alpha}_i\ga'_ix_0)_{i\in\NN}$ where
$(\beta_i)_{i\in\NN}$, $(\overline{\alpha}_i)_{i\in\NN}$,
$(\ga'_i)_{i\in\NN}$ are sequences in respectively $\Ga_{p,i}$,
$\Ga_{p,i} \bs\!\bs\Ga_p$ and $\Ga_{p}\bs\!\bs\Ga$. Up to extraction,
if $\xi\neq p$, since the limit set of $\Ga_p$ is reduced to $\{p\}$,
we may assume that $\lim_{i\ra+\infty}\ga'_ix_0=\xi'\in\partial_\infty
\wt M$. By Lemma \ref{lem:quasigeod} (i), the points $\ga'_ix_0$
belong to the union of the geodesic lines starting from $p$ and
passing through the closed ball $\overline{B}(x_0,c_1)$. The set of
endpoints of these geodesic lines is closed and does not contain
$p$. Hence $\xi'\neq p$. Since any compact neighbourhood of $\xi'$ not
containing $p$ is mapped into any given neighbourhood of $p$ by all
except finitely many elements of $\Ga_p$, if the sequence $(\beta_i
\overline{\alpha}_i)_{i\in\NN}$ in $\Ga_p$ takes infinitely many
values, then $\xi=p$. Hence up to extraction, if $\xi\neq p$, the
sequence $(\beta_i \overline{\alpha}_i)_{i\in\NN}$ is constant, and so
is $(\beta_i)_{i\in\NN}$: therefore $\xi\in \beta_0 \F^+_{p,i}$. This
proves the result. \cqfd

\medskip Let $\F_{p,i}=\nu P_{D_i}(\F^+_{p,i}-\partial_\infty D_i)\cap
\pi^{-1}(\H_p)$.  The images of $\F_{p,i}$ under the elements of
$\Ga_{p,i}$ cover $\pi^{-1}(\H_p)\cap \operatorname{Supp}
\wt\sigma_{D_i}$. It follows from Lemma \ref{lem:quasigeod} (ii) that
there exists $c_3>0$ such that $P_{D_i}(\F_{p,i}^+)=\pi(\F_{p,i})$ is
contained in the $c_3$-neighbourhood of the geodesic ray $[x_0,p[\,$.

In order to prove the finiteness of $\pi_*\sigma_{\D}(\Ga_p\bs\H_p)$,
we hence only have to prove the finiteness of $\wt\sigma_{D_i}
(\F_{p,i})$.

\medskip The next lemma, which uses the assumption that $\Ga$ is of
divergence type, gives a control on the Patterson measure $\mu_y$ of
$\F_{p,i}^+$ as $y$  converges radially to $p$.

\blemm \label{lem:decaypatterson}
There exists $c_4>0$ such that for every $y$ on the geodesic
ray $[x_0,p[\,$, we have
$$
\mu_y(\F_{p,i}^+)
\leq c_4\;e^{(2(\delta_{p}-\delta_{p,i})-\delta)\,d(x_0,\,y)}\;.
$$
\elemm

\dem For all $s\geq 0$ and $y\in\wt M$, for every subgroup $\Ga'$ of
$\Ga$, let
$$
P_{\Ga',\,y}(s)=\sum_{\ga\in\Ga'}\;e^{-s\,d(y,\,\ga x_0)}\in[0,+\infty]\;,
$$
and let $\D_y$ be the unit Dirac mass at the point $y$.  Since
$\Ga$ is of divergence type, the Patterson measure $\mu_y$ is the
weak-star limit as $s\ra\delta^+$ of the measures
$$
\mu_{y,s}=\frac{1}{P_{\Ga,\,x_0}(s)}
\sum_{\ga\in\Ga}\;e^{-s\,d(y,\,\ga x_0)} \D_{\ga x_0}
$$ 
(see for instance \cite{Roblin03}). By discreteness and Lemma
\ref{lem:eq:decomplimset}, there exists a finite subset $F$ of
$\Ga_{p,i}$ such that $\bigcup_{\beta\in F} \beta\F_{p,i}^+-\{p\}$ is
a neighbourhood of $\F_{p,i}^+-\{p\}$ in $\Lambda\Ga-\{p\}$. Since
$\Ga$ is of divergence type, the measure $\mu_y$ has no atom at $p$
(see for instance \cite[Coro.~1.8]{Roblin03}).  Hence there exists
$c>0$ such that for every $y\in [x_0,p[\,$, we have
$$
\mu_y(\F_{p,i}^+)\leq c\lim_{s\ra \delta^+}\frac{1}{P_{\Ga,\,x_0}(s)}
\sum_{\ov\alpha\in \Ga_{p,i}\bs\!\bs\Ga_p,\;\ga'\in\Ga_{p}\bs\!\bs\Ga}
\;e^{-s\,d(y,\,\ov\alpha\ga' x_0)}\;.
$$

Let 
$$
Q_y(s)= \sum_{\ov\alpha\in \Ga_{p,i}\bs\!\bs\Ga_p}
\;e^{-s\,d(y,\,\ov\alpha y)}\;\;{\rm and}\;\;
R(s)=\sum_{\ga'\in\Ga_{p}\bs\!\bs\Ga}
\;e^{-s\,d(x_0,\,\ga' x_0)}\;.
$$
By the lower bound in Lemma \ref{lem:quasigeod}  (i), we have
$$
\sum_{\ov\alpha\in \Ga_{p,i}\bs\!\bs\Ga_p,\;\ga'\in\Ga_{p}\bs\!\bs\Ga}
\;e^{-s\,d(y,\,\ov\alpha\ga' x_0)}\leq
e^{s\,c_1}\;e^{-s\,d(y,\,x_0)}\;Q_y(s)\;R(s)\;.
$$
Similarly, by the upper bound in Lemma \ref{lem:quasigeod} (i), we
have $P_{\Ga,\,x_0}(s)\geq P_{\Ga_p,\,x_0}(s)R(s)$. We will prove below
that the series $Q_y(\delta)$ converges. Thus, 
even if $P_{\Ga_p,\,x_0}(\delta)=+\infty$, we have
\begin{equation}\label{eq:majomesPatteny}
\mu_y(\F_{p,i}^+)\leq
\frac{c\;e^{\delta\,c_1}}{P_{\Ga_p,\,x_0}(\delta)}
\;e^{-\delta\,d(y,\,x_0)}\;Q_y(\delta)\;.
\end{equation}

By the convexity of the horoball of center $p$ whose boundary contains
$y$ and by standard quasi-geodesic arguments, there exist two
constants $c',c''>0$ such that for every $\ov\alpha\in\Ga_{p,i} \bs\!
\bs\Ga_p$, if $d(y,\ov\alpha y)> c'$ then
$$
d(y,\ov\alpha y)+2d(y,x_0)- c''\leq d(x_0,\ov\alpha x_0) 
\leq d(y,\ov\alpha y)+2d(y,x_0)\;.
$$ 
If $d(y,\ov\alpha y)\leq c'$ then $d(x_0,\ov\alpha x_0) \leq
2d(y,x_0)+c'$, by the triangle inequality.  We hence have, using the
notation $t\mapsto h(t)$ introduced in the proof of Lemma
\ref{lem:expcritrelatif},
\begin{align*}
Q_y(\delta)&=\sum_{\tiny
\begin{array}{c} \ov\alpha\in \Ga_{p,i}\bs\!\bs\Ga_p\\
d(y,\,\ov\alpha y)\leq c'\end{array}}
\;e^{-\delta\,d(y,\,\ov\alpha y)}+ \sum_{\tiny
\begin{array}{c} \ov\alpha\in \Ga_{p,i}\bs\!\bs\Ga_p\\
d(y,\,\ov\alpha y)> c'\end{array}}
\;e^{-\delta\,d(y,\,\ov\alpha y)}\\&\leq
\sum_{\tiny
\begin{array}{c} \ov\alpha\in \Ga_{p,i}\bs\!\bs\Ga_p\\
d(x_0,\,\ov\alpha x_0)\leq 2\,d(y,\,x_0)+c' \end{array}}
\;1+\sum_{\tiny
\begin{array}{c} \ov\alpha\in \Ga_{p,i}\bs\!\bs\Ga_p\\
d(x_0,\ov\alpha x_0) \geq 2\,d(y,x_0)+c'-c''\end{array}}
\;e^{-\delta\,d(x_0,\,\ov\alpha x_0)+2\,\delta \,d(y,\,x_0)}
\\&\leq h\big(2\,d(y,\,x_0)+c'\big)\;+\;e^{2\,\delta\, d(y,\,x_0)}
\sum_{n\geq 2\,d(y,\,x_0)+c'-c''-1}\;h(n+1)\;e^{-\delta n}\;.
\end{align*}
Since $h(t)\leq c_2\,e^{(\delta_p-\delta_{p,i})t}$ by Lemma
\ref{lem:expcritrelatif}, and by a geometric series summation argument
since $\delta-\delta_p+\delta_{p,i}>0$ by the assumption
\eqref{eq:condiexpocrit}, we therefore have
\begin{align*}
Q_y(\delta)&\leq 
c_2\,e^{(\delta_p-\delta_{p,i})c'}e^{2(\delta_p-\delta_{p,i})d(y,\,x_0)}
\\&\;\;\;\;+c_2\,e^{\delta_p-\delta_{p,i}+2\,\delta\,d(y,\,x_0)}
\sum_{n\geq 2\,d(y,\,x_0)+c'-c''-1}e^{(\delta_p-\delta_{p,i}-\delta) n}
\\&\leq c'''e^{2(\delta_p-\delta_{p,i})d(y,\,x_0)}\;,
\end{align*}
for some $c'''>0$. Using Equation \eqref{eq:majomesPatteny}, this
proves Lemma \ref{lem:decaypatterson}.  
\cqfd

\medskip Let $\rho:[0,+\infty[\;\ra\wt M$ be the geodesic ray with
origin $x_0$ and point at infinity $p$. For every $n\in\NN$, let $A_n$
be the set of points $\xi\in\partial_\infty\wt M-\{p\}$ such that the
closest point to $\xi$ on the geodesic ray $\rho$ belongs to
$\rho([n,n+1])$. Note that $\bigcup_{n\in\NN}A_n=\partial_\infty\wt
M-\{p\}$.

\blemm \label{lem:majoshadow}
There exists $c_5>0$ such that for every $n\in\NN$, we have
$$
\mu_{x_0}(\F_{p,i}^+\cap A_n)
\leq c_5\;e^{2(\delta_{p}-\delta_{p,i}-\delta)n}\;.
$$
\elemm

\dem For every $\xi\in A_n$, since the angle at $\rho(n)$ between
$[\rho(n),x_0]$ and $[\rho(n),\xi[$ is at least $\frac{\pi}{2}$ if
$n\neq 0$, we have $\beta_\xi(x_0,\rho(n))\geq d(x_0,\rho(n))-
2\varepsilon= n-2\varepsilon$.  Hence, by Equation \eqref{eq:RNdensity}
and by Lemma \ref{lem:decaypatterson}, we have
\begin{align*}
\mu_{x_0}(\F_{p,i}^+\cap A_n)&=\int_{\xi\in \F_{p,i}^+\cap A_n}
e^{-\delta\,\beta_\xi(x_0,\,\rho(n))}\;d\mu_{\rho(n)}(\xi)
\leq e^{-\delta n +2\delta\varepsilon}\;\mu_{\rho(n)}(\F_{p,i}^+)\\&
\leq c_4\,e^{2\delta\varepsilon}\,e^{2(\delta_{p}-\delta_{p,i}-\delta)n}\;.\;\Box
\end{align*}

\medskip After this series of lemmas, let us prove the finiteness of
$\wt\sigma_{D_i}(\F_{p,i})$, which concludes the proof of Theorem
\ref{theo:geneohshah}.

With $a_n= \wt\sigma_{D_i}(\F_{p,i}\cap \nu P_{D_i}(A_n))$, we only
have to prove that the series $\sum_{n\in\NN} a_n$ converges. For
every $\xi\in \F_{p,i}^+\cap A_n$, by the definition of $c_3$, the
point $P_{D_i}(\xi)$ lies at distance less than a constant from
$\rho([n,n+1])$. Hence there exists a constant $c>0$ such that
$\beta_\xi(P_{D_i}(\xi),x_0)\geq -n-c$. By the definition of the
skinning measures in Equation \eqref{eq:defiskinmeas} and by Lemma
\ref{lem:majoshadow}, we hence have
\begin{align*}
\wt\sigma_{D_i}(\F_{p,i}\cap \nu P_{D_i}(A_n))&\leq
\int_{\xi\in \F_{p,i}^+\cap A_n}
e^{-\delta\,\beta_\xi(P_{D_i}(\xi),\,x_0)}\;d\mu_{x_0}(\xi)
\leq e^{\delta n +\delta c}\;\mu_{x_0}(\F_{p,i}^+\cap A_n)\\&
\leq c_5\,e^{\delta c}\,e^{(2(\delta_{p}-\delta_{p,i})-\delta)n}\;.
\end{align*}
By the assumption $\delta>2(\delta_{p}-\delta_{p,i})$ in Equation
\eqref{eq:condiexpocrit}, a geometric series summation argument
proves that $\sum_{n\in\NN} a_n$ converges.  This completes the proof
of Theorem \ref{theo:geneohshah}.  
\cqfd

\section{Equidistribution of equidistant submanifolds}
\label{sec:equidistribution}

Let $\wt M,\Ga,x_0,M$ and $T^1M$ be as in Section
\ref{sec:geometry}. Assume that the critical exponent $\delta_\Ga$ of
$\Ga$ is finite. Let $(\mu_x)_{x\in\wt M}$ be a Patterson density of
dimension $\delta_\Ga$, with Bowen-Margulis measures $\wt m_{\rm BM}$
and $m_{\rm BM}$ on $T^1\wt M$ and $T^1M$, respectively.  Let
$\C=(C_i)_{i\in I}$ be a $\Ga$-equivariant family of proper nonempty
closed convex subsets of $\wt M$.  Let $\C_{t}=(\N_tC_i)_{i\in I}$ (in
particular $\C_0=\C$), and let $\wt \sigma_{\C_t}$ and $\sigma_{\C_t}$
be the skinning measures of $\C_t$ on $T^1\wt M$ and $T^1M$,
respectively.  Let $\Omega=(\Omega_i)_{i\in I}$ be a locally finite
$\Ga$-equivariant family of subsets of $T^1\wt M$, where $\Omega_i$ is a
measurable subset of $\normal C_i$ with $\wt\sigma_{C_i}
(\partial\Omega_i) =0$ for every $i\in I$. Let $\sim\;=\;
\sim_{\Omega}$ be the equivalence relation on $I$ defined at the end
of Section \ref{sec:skinning}.  As we have already defined when
$\Omega=\C$, let
$$
\wt\sigma_\Omega=\sum_{i\in I/\sim}\wt\sigma_{C_i}|_{\Omega_i}\;,
$$ 
which is a $\Ga$-invariant locally finite positive Borel measure on
$T^1\wt M$ (independent of the choice of representatives in $I/\sim$).
Hence, $\wt\sigma_\Omega$ induces a locally finite positive Borel
measure $\sigma_\Omega$ on $T^1M$.  Note that $\flow{t}\Omega_i\subset
\normal\N_t C_i$ and as in the end of Section \ref{sec:skinning}, for
every $t>0$, we have that
\begin{equation}\label{eq:poussskinflot}
\|\sigma_{\flow{t}\Omega}\|=e^{\delta_\Ga t}\|\sigma_{\Omega}\|\,.
\end{equation}

The aim of this section is to prove, under some finiteness
assumptions, that the measures $\sigma_{\flow{t}\Omega}$ on $T^1M$
equidistribute to the Bowen-Margulis measure on $T^1M$ as 
$t\ra+\infty$. We start by introducing the test functions
approximating the support of the measures  $\sigma_{\flow{t}\Omega}$.

\medskip Assume that the number of orbits of $\Ga$ on the set of
elements $i\in I$, such that the intersection $\partial_\infty C_i\cap
\Lambda\Ga$ is empty, is finite (this condition is stronger than the
requirement on $\C$ to be locally finite).
Under this assumption, by Lemma \ref{lem:veryoldlemseven}, there
exists $R>0$ such that for every $i\in I$, for every $w\in \normal
C_i$, we have $\mu^{\rm ss}_w(V_{w,R})>0$, where $V_{w,R}$ is the open
ball of radius $R$ and center $w$ for the Hamenst\"adt distance in the
strong stable leaf $W^{\rm ss}(w)$. We fix such an $R$.

For every $\eta>0$, let $h_{\eta,\,R}: T^1\wt M\ra[0,+\infty]$ be
the measurable map defined by
$$
h_{\eta,\,R}(w)=\frac{1}{2\eta\;\mu^{\rm ss}_w(V_{w,\,R})}\;.
$$
Note that $h_{\eta,\,R}$ is $\Ga$-invariant by Equation
\eqref{eq:invarcondstab} and that $h_{\eta,\,R} \circ
\flow{-t}=e^{-\delta_\Ga t}\;h_{\eta,\,e^{-t}R}$ for every $t\in\RR$:
indeed, for every $w\in T^1\wt M$, we have by Equation
\eqref{eq:invargeoflocondstab}
$$
h_{\eta,\,R}(\flow{-t}w)=
\frac{1}{2\eta\;\mu^{\rm ss}_{\flow{-t}w}(V_{\flow{-t}w,\,R})}
=\frac{1}{2\eta\;e^{\delta_\Ga t}
\mu^{\rm ss}_{w}(\flow{t}V_{\flow{-t}w,\,R})}
=\frac{e^{-\delta_\Ga t}}{2\eta\;\mu^{\rm ss}_{w}(V_{w,\,e^{-t}R})}
\;.
$$
For every $i\in I$, let $\V_{\eta,\,R,\,i}= \V_{\eta,\,R}(\Omega_i)$ be the
dynamical thickening of $\Omega_i$ defined at the end of Section
\ref{sec:geometry}. Note that $\ga\V_{\eta,\,R,\,i}=\V_{\eta,\,R,\,\ga i}$ for
every $\ga\in \Ga$ and every $i\in I$. 

We denote by $\chi_A$ the characteristic function of a subset $A$.  We
will use the test function $\wt\phi_\eta= \wt\phi_{\eta,\,R,\,\Omega}:
T^1\wt M\to[0,+\infty[$ defined by (using the convention $\infty\times
0=0$)
$$
\wt\phi_\eta(v)=\sum_{i\in I/\sim}
h_{\eta,\,R}\circ f_{C_i}(v)\;\chi_{\V_{\eta,\,R,\,i}}(v)\;,
$$  
where $f_{C_i}:U_{C_i}\to\normal C_i$ is the fibration defined in
Section \ref{sec:geometry}. Note that $\V_{\eta,\,R,\,i}$ is contained
in $U_{C_i}$, and we define $h_{\eta,\,R}\circ f_{C_i}(v)\;
\chi_{\V_{\eta,\,R,\,i}}(v) =0$ if $v\notin\V_{\eta,\,R,\,i}$.

\blemm The function $\wt\phi_\eta$ is well defined, measurable and
$\Ga$-invariant. Furthermore, for every $t\in[0,+\infty[\,$, we have
$\wt\phi_{\eta,\,R,\,\Omega}\circ \flow{-t}=e^{-\delta_\Ga t}\;
\wt\phi_{\eta,\,e^{-t}R,\,\flow{t}\Omega}$.  
\elemm

\dem The function $\wt\phi_\eta$ is well defined, since
${\Omega_i}= {\Omega_j}$ and $\V_{\eta,\,R,\,i}= \V_{\eta,\,R,\,j}$ if
$i\sim j$, since $h_{\eta,\,R}\circ f_{C_i}(v)$ is finite if $v\in
\V_{\eta,\,R,\,i}$ (by the definition of $R$), and since the sum
defining $\wt\phi_\eta(v)$ has only finitely many nonzero terms, by the
local finiteness of the family $\Omega$ (given $v$, the summation over
$I/\!\!\sim$ may be replaced by a summation over the finite set
$\{i\in I: v\in \V_{\eta,\,R,\,i}\}/\!\!\sim$).

The function $\wt\phi_\eta$ is $\Ga$-invariant since
$\chi_{\V_{\eta,\,R,\,i}}\circ \ga= \chi_{\ga^{-1}\V_{\eta,\,R,\,i}}=
\chi_{\V_{\eta,\,R,\,\ga^{-1}i}}$ and 
$$
h_{\eta,\,R} \circ f_{C_i} \circ \ga=
h_{\eta,\,R} \circ \ga \circ f_{\ga^{-1}C_i} =
h_{\eta,\,R}\circ f_{C_{\ga^{-1}i}}
$$
and by a change of index in the above sum. 

Let $t\geq 0$. The last claim follows by noting that 
 $$ 
\chi_{\V_{\eta,\,R}(\Omega_i)}\circ \flow {-t}=
\chi_{\flow {t}\V_{\eta,\,R}(\Omega_i)}
 =\chi_{\V_{\eta,\,e^{-t}R}(\flow{t}\Omega_i)}\;,
$$
and
$$
h_{\eta,\,R}\circ f_{C_i}\circ \flow {-t}=h_{\eta,\,R}\circ f_{C_i} =
h_{\eta,\,R}\circ\flow {-t}\circ \flow {t} \circ f_{C_i}=
e^{-\delta_\Ga t}h_{\eta,\,e^{-t}R}\circ f_{\N_t C_i}\;. \;\;\;\Box
$$

\medskip Hence the test function $\wt\phi_\eta$ defines, by passing to
the quotient, a measurable function $\phi_\eta=\phi_{\eta,\,R,\,\Omega}
: T^1 M\to[0,+\infty[\,$, such that for every $t\in[0,+\infty[\,$,
we have
\begin{equation}\label{eq:actflotestbas}
\phi_{\eta,\,R,\,\Omega}\circ \flow{-t}=
e^{-\delta_\Ga t}\;\phi_{\eta,\,e^{-t}R,\,\flow{t}\Omega}\;.
\end{equation}

\bprop \label{prop:integrable} Assume that the Bowen-Margulis measure
of $T^1M$ is finite. For every $\eta>0$, we have
$\int\phi_\eta\,dm_{\rm BM}=\|\sigma_{\Omega}\|$.  In particular, the
function $\phi_\eta$ is integrable for the Bowen-Margulis measure if
and only if the measure $\sigma_\Omega$ is finite.  \eprop

\dem Let $i\in I$ and let $K_i$ be a measurable subset of
$\Omega_i$. By the disintegration result of Proposition 
\ref{prop:disintegration} (more precisely by Equation
\eqref{eq:decomskinstrongstab}), and by the definitions of the
function $h_{\eta,\,R}$ and of the set $\V_{\eta,\,R,\,i}=
\bigcup_{w\in\Omega_i} \bigcup_{s\in\;]-\eta,\eta[} g^sV_{w,\,R}$,
we have
\begin{align*}
\int_{\V_{\eta,\,R,\,i}\cap f^{-1}_{C_i}(K_i)}h_{\eta,\,R}\circ f_{C_i}\; 
d\wt m_{\rm BM}&=
\int_{w\in K_i} h_{\eta,\,R}(w)\int_{v'\in V_{w,R}} \int_{-\eta}^\eta
\;ds\,d\muss{w}(v')\,d\wt\sigma_{C_i}(w)\\ &=\wt\sigma_{C_i}(K_i)\;. 
\end{align*}
Let $\Delta_\Ga$ be a {\it fundamental domain} for the action of $\Ga$
on $T^1\wt M$, that is, $\Delta_\Ga$ is the closure of its interior,
its boundary has measure $0$ for the Bowen-Margulis measure, the
images of $\Delta_\Ga$ by the elements of $\Ga$ have pairwise disjoint
interiors and cover $T^1\wt M$, and any compact subset of $T^1\wt M$
meets only finitely many images of $\Delta_\Ga$ by elements of
$\Ga$. Such a fundamental domain exists since the Bowen-Margulis
measure of $T^1M$ is finite (see 
for instance \cite[page 13]{Roblin03}). By the definition of the test
function $\wt\phi_\eta$, we have
$$
\int_{T^1M}\phi_\eta\;dm_{\rm BM}=
\int_{\Delta_\Ga}\wt\phi_\eta\;d\wt m_{\rm BM}=
\sum_{i\in I/\sim} \int_{\V_{\eta,\,R,i}\cap\Delta_\Ga} 
h_{\eta,\,R}\circ f_{C_i}\;d\wt m_{\rm BM}\;.
$$
By the definition of the  measure $\sigma_\Omega$, we have
$$
\|\sigma_\Omega\|= \wt\sigma_\Omega(\Delta_\Ga)=
\sum_{i\in I/\sim} \wt\sigma_{C_i}(\Delta_\Ga\cap \Omega_i)\;.
$$
By an easy multiplicity argument, the result follows. \cqfd

\medskip Now, we can state and prove the main result of this paper.

\btheo\label{theo:equid} Let $\wt M$ be a complete simply connected
Riemannian manifold with sectional curvature bounded above by
$-1$. Let $\Ga$ be a discrete, nonelementary group of isometries of
$\wt M$, with finite critical exponent. Assume that the Bowen-Margulis
measure $m_{\rm BM}$ of $\Ga$ on $T^1M$ is finite and mixing for the
geodesic flow.  Let $\C=(C_i)_{i\in I}$ be a $\Ga$-equivariant family
of nonempty proper closed convex subsets of $\wt M$. Let $\Omega=
(\Omega_i)_{i\in I}$ be a locally finite $\Ga$-equivariant family of
measurable subsets $\Omega_i\subset\normal C_i$ with
$\wt\sigma_{C_i}(\partial\Omega_i)=0$.  Assume that $\sigma_\Omega$ is
finite and nonzero. Then, as $t\ra+\infty$,
$$
\frac{1}{\|\sigma_{\flow{t}\Omega}\|}\,\sigma_{\flow{t}\Omega}\;\;
\stackrel{*}{\rightharpoonup}
\;\;\frac{1}{\|m_{\rm BM}\|}\,m_{\rm BM}\;. 
$$
\etheo

In particular, if $\C=(C_i)_{i\in I}$ is a locally finite
$\Ga$-equivariant family of nonempty proper closed convex subsets of
$\wt M$ with finite nonzero skinning measure, then the skinning
measure $\sigma_{\C_t}$ on $T^1M$ of $\C_{t}=(\N_tC_i)_{i\in I}$
equidistributes to the Bowen-Margulis measure 
as
$t\ra+\infty$.

\medskip
\dem Given three numbers $a,b,c$ (depending on some parameters), we
write $a=b \pm c$ if $|a-b|\leq c$.

We may assume that $\Ga\bs I$ is finite. Indeed, if $\ov J$
is a big enough finite subset of $\Ga\bs I$, if $J$ is
the preimage of $\ov J$ by the canonical map $I\ra \Ga\bs I$, since
the  measure  $\sigma_\Omega$ is finite, the contribution of the family
$(\flow{t}\Omega_i)_{i\in I-J}$ is negligible compared to the one of
$(\flow{t}\Omega_i)_{i\in J}$ (they grow at equal rate as $t$ tends to
$+\infty$, by Equation \eqref{eq:poussskinflot}). 

Hence we may consider $R>0$ as was fixed in the beginning of Section
\ref{sec:equidistribution} and, for every $\eta>0$, the 
test function $\phi_\eta=\phi_{\eta,\,R,\,\Omega}$ as defined above.

Fix $\psi\in\C_c(T^1M)$. Let us prove that 
$$
\lim_{t\ra+\infty}\;\frac{1}{\|\sigma_{\flow{t}\Omega}\|}\;
\int_{T^1M}\psi\;d\sigma_{\flow{t}\Omega}=\frac{1}{\|m_{\rm BM}\|}\;
\int_{T^1M}\psi\;dm_{\rm BM}\;.
$$
Given a fundamental domain $\Delta_\Ga$ for the action of $\Ga$ on
$T^1\wt M$ as above, by a standard argument of finite partition of
unity, we may assume that there exists a map $\wt\psi: T^1\wt
M\ra\RR$ whose support is contained in $\Delta_\Ga$ such that $\wt
\psi=\psi\circ p$, where $p:T^1\wt M\ra\Ga\bs T^1\wt M$ is the
canonical projection (which is $1$-Lipschitz). Fix $\epsilon>0$. Since
$\wt\psi$ is uniformly continuous, for every $\eta>0$ small enough,
and for every $t\geq 0$ big enough, for every $w\in T^1\wt M$ and
$v\in V_{w,\eta,e^{-t}R}$, we have
\begin{equation}\label{eq:unifcont}
\wt\psi(v)= \wt\psi(w)\pm \frac{\epsilon}{2}\;.
\end{equation}

We have, using respectively 
\begin{itemize}
\item Equation \eqref{eq:actflotestbas} and the definition of $\wt
  \psi$ for the first and second equality,
\item the definition of the test function $\wt \phi_\eta$ for the
third equality,
\item Equation \eqref{eq:unifcont} and the fact that the support
of $\wt \psi$ is contained in $\Delta_\Ga$ for the fourth
equality,
\item the invariance of the Bowen-Margulis measure under the geodesic
  flow, and Equation \eqref{eq:decomskinstrongstab} as in the proof of
  Proposition \ref{prop:integrable} for the fifth equality,
\item the definition of $h_{\eta,\,e^{-t}R}$ and Proposition
  \ref{prop:integrable} for the sixth equality:
\end{itemize}

\begin{align*}
&\int_{T^1M} \phi_\eta\circ\flow {-t}\;\psi\;dm_{\rm BM}\\
=\;&
e^{-\delta_\Ga t}\int_{T^1M} \phi_{\eta,\,e^{-t}R,\flow{t}\Omega}\;\psi\;dm_{\rm BM}=
e^{-\delta_\Ga t}\int_{T^1\wt M} \wt \phi_{\eta,\,e^{-t}R,\flow{t}\Omega}\;\wt\psi\;
d\wt m_{\rm BM}\\
=\;&e^{-\delta_\Ga t}\sum_{i\in I/\sim}\;
\int_{\V_{\eta,\,e^{-t}R}(\flow{t}\Omega_i)}\; 
h_{\eta,\,e^{-t}R}\circ f_{\N_tC_i}\;\wt\psi\;d\wt m_{\rm BM}\\
=\;& 
e^{-\delta_\Ga t}\sum_{i\in I/\sim}\;
\int_{\V_{\eta,\,e^{-t}R}(\flow{t}\Omega_i)}\; 
(h_{\eta,\,e^{-t}R}\;\wt\psi)\circ f_{\N_tC_i}\;d\wt m_{\rm BM} 
\pm\frac{\epsilon}{2}\int_{\Delta_\Ga} \wt\phi_\eta\circ\flow t\;
d\wt m_{\rm BM}\\
=\;& e^{-\delta_\Ga t}\sum_{i\in I/\sim}\;
\int_{w\in \flow{t}\Omega_i}\; 
h_{\eta,\,e^{-t}R}(w)\;\wt\psi(w)\;(2\eta)\;\mu^{\rm ss}_w(V_{w,\,e^{-t}R})
\;d\wt\sigma_{\N_tC_i} \pm\frac{\epsilon}{2}\int_{T^1M} 
\phi_\eta\;dm_{\rm BM} \\
=\;&
e^{-\delta_\Ga t}\sum_{i\in I/\sim}\;
\int_{w\in \flow{t}\Omega_i}\; \wt\psi(w)
\;d\wt\sigma_{\N_tC_i} \pm\frac{\epsilon}{2}\,\|\sigma_\Omega\|\\=
\;&
e^{-\delta_\Ga t}\int \psi\;d\sigma_{\flow{t}\Omega} 
\pm\frac{\epsilon}{2}\,\|\sigma_\Omega\|\;.
\end{align*}
Hence, using Equation \eqref{eq:poussskinflot} for the first equality,
the previous computation for the second equality, the invariance of
the Bowen-Margulis measure under the geodesic flow for the third
equality, and Proposition \ref{prop:integrable} for the last one, we
have, for $\eta>0$ small enough and $t\geq 0$ big enough,
\begin{align} 
  \frac{\int
    \psi\;d\sigma_{\flow{t}\Omega}}{\|\sigma_{\flow{t}\Omega}\|}& 
= \frac{\int
    \psi\;d\sigma_{\flow{t}\Omega}}{e^{\delta_\Ga t}\|\sigma_{\Omega}\|}=
  \frac{\int_{T^1M} \phi_\eta\circ\flow {-t}\;\psi\;dm_{\rm
      BM}}{\|\sigma_{\Omega}\|} \pm\frac{\epsilon}{2}\nonumber\\ &
= \frac{\int_{T^1M}
    \phi_\eta\;\psi\circ\flow t\;dm_{\rm BM}}{\|\sigma_{\Omega}\|}
  \pm\frac{\epsilon}{2}= 
\frac{\int_{T^1M} \phi_\eta\;\psi\circ\flow t\;dm_{\rm
      BM}}{\int_{T^1M}\phi_\eta\;dm_{\rm BM}} \pm\frac{\epsilon}{2} \;.
\label{eq:prepamixing}
\end{align}
By the mixing property of the geodesic flow on $T^1M$, for $t\geq 0$
big enough (while $\eta$ is fixed), we have 
$$
\frac{\int_{T^1M}
  \phi_\eta\;\psi\circ\flow t\;dm_{\rm
    BM}}{\int_{T^1M}\phi_\eta\;dm_{\rm BM}}= \frac{\int_{T^1M}
  \psi\;dm_{\rm BM}}{\|m_{\rm BM}\|}\pm\frac{\epsilon}{2}\;.
$$  
This proves the result.
\cqfd

\medskip We conclude this section by proving Theorem \ref{theo:intro}
in the introduction. The definition of a properly immersed locally
convex subset is recalled in the beginning of the proof.

\medskip\noindent{\bf Proof of Theorem \ref{theo:intro}. } Let $M, C$
be as in the statement of Theorem \ref{theo:intro}, that is, they satisfy
the following property. Let $\wt M\ra M$ be a universal covering of
$M$, with covering group $\Ga$. Let $\wt C\ra C$ be a covering map
which is a universal covering over each component of $C$. The
immersion from $C$ to $M$ lifts to an immersion $f: \wt C\ra\wt M$,
which is, on each connected component of $\wt C$, an embedding whose
image is a convex subset of $\wt M$.

Let $I=\Ga\times \pi_0\wt C$ with the action of $\Ga$ defined by
$\ga\cdot (\alpha,c) = (\ga \alpha,c)$ for all $\ga,\alpha\in \Ga$ and
every component $c$ of $\wt C$. Consider the family $\C=(C_i)_{i\in
  I}$ where $C_i= \alpha \;\overline{f(c)}$ if $i=(\alpha,c)$.  Then
$\C$ is a $\Ga$-equivariant family of nonempty closed convex subsets
of $\wt M$, which is locally finite since $C$ is properly immersed in
$M$. The result then follows from Theorem \ref{theo:equid}.  
\cqfd

\section{Exponential rate of equidistribution}
\label{sec:expo}

Let $\wt M,\Ga,M,T^1M,m_{\rm BM},\C,\C_{t}$ and $\sigma_{\C_t}$ be as
in the beginning of Section \ref{sec:equidistribution}. When the
Bowen-Margulis measure $m_{\rm BM}$ is finite, we denote by
$\overline{m}_{\rm BM}$ its normalisation to a probability measure.

In this section, we show, under the finiteness assumptions of Theorem
\ref{theo:equid}, that in the known cases when the geodesic flow is
exponentially mixing, the skinning measure equidistributes to the
Bowen-Margulis measure with exponential speed.  To begin with, we
recall the two types of exponential mixing results that are available.
In order to prove our estimates for the rate of equidistribution using
these results, we will smoothen (accordingly to the two regularities)
our test function $\phi_\eta$ defined in the previous section.

\medskip Firstly, when $M$ is locally symmetric with finite volume,
then the boundary at infinity of $\wt M$, the strong unstable,
unstable, stable, and strong stable foliations of $T^1\wt M$ are
smooth. Hence, for all $\ell\in\NN$, talking about $\C^\ell$-smooth
leafwise defined functions on $T^1M$ makes sense. We will denote by
$\C_c^\ell(T^1M)$ the vector space of $\C^\ell$ smooth functions on
$T^1M$ with compact support and by $\|\psi\|_\ell$ the Sobolev
$W^{\ell,2}$-norm of any $\psi\in\C_c^\ell(T^1M)$. Note that now the
Bowen-Margulis measure $m_{\rm BM}$ of $T^1M$ is the unique (up to a
multiplicative constant) locally homogeneous smooth measure on $T^1M$
(hence it coincides, up to a multiplicative constant, with the Liouville
measure).

Given $\ell\in \NN$, we will say that the geodesic flow on $T^1M$ is
{\it exponentially mixing for the Sobolev regularity $\ell$} (or that
it has {\em exponential decay of $\ell$-Sobolev correlations}) if
there exist $c,\kappa>0$ such that for all $\phi,\psi\in
\C_c^\ell(T^1M)$ and all $t\in\RR$, we have
$$
\Big|\int_{T^1M} \phi\circ g^{-t}\;\psi\;
d\overline{m}_{\rm BM}-\int_{T^1M}
\phi\;d\overline{m}_{\rm BM}\int_{T^1M} \psi\;
d\overline{m}_{\rm BM}\;\Big|\leq
c\,e^{-\kappa |t|}\;\|\psi\|_\ell\;\|\phi\|_\ell\;.
$$
When $\Ga$ is a arithmetic lattice in the isometry group of $\wt M$,
this property, for some $\ell\in\NN$, follows from
\cite[Theorem~2.4.5]{KleMar96}, with the help of \cite[Theorem
3.1]{Clozel03} to check its spectral gap property, and of
\cite[Lemma~3.1]{KleMar99} to deal with finite cover problems.

Secondly, when $\wt M$ has pinched negative sectional curvature with
bounded derivatives, then the boundary at infinity of $\wt M$, the
strong unstable, unstable, stable, and strong stable foliations of
$T^1\wt M$ are only H\"older-smooth (see for instance \cite{Brin95}
when $\wt M$ has a compact quotient, and
\cite[Theo.~7.3]{PauPolSha11}).  Hence the appropriate regularity on
functions on $T^1\wt M$ is the H\"older one.  For every
$\alpha\in\;]0,1[\,$, we denote by $\operatorname{C}_{\rm c} ^\alpha
(X)$ the space of $\alpha$-H\"older continuous real-valued functions
with compact support on a metric space $(X,d)$, endowed with the
H\"older norm
$$
\|f\|_\alpha=
\|f\|_\infty+\sup_{x,\,y\in X,\;x\neq y}\frac{|f(x)-f(y)|}{d(x,y)^\alpha}\,.
$$

Assuming the Bowen-Margulis measure $m_{\rm BM}$ on $T^1M$ to be
finite, given $\alpha\in\;]0,1[$, we will say that the geodesic flow
on $T^1M$ is {\em exponentially mixing for the H\"older regularity
  $\alpha$} (or that it has {\em exponential decay of
  $\alpha$-H\"older correlations}) if there exist $\kappa,c >0$ such
that for all $\phi,\psi\in \operatorname{C}_{\rm c}^\alpha(T^1M)$ and
all $t\in\RR$, we have
$$
\Big|\int_{T^1M}\phi\circ\flow{-t}\;\psi\;d\overline{m}_{\rm BM}-
\int_{T^1M}\phi\; d\overline{m}_{\rm BM} 
\int_{T^1M}\psi\;d\overline{m}_{\rm BM}\;\Big|
\le c\;e^{-\kappa|t|}\;\|\phi\|_\alpha\;\|\psi\|_\alpha\,.
$$
This holds if $ M$ is compact and has dimension $2$ by the work of Dolgopyat
\cite{Dolgopyat98} or if $M$ is compact and  locally symmetric by
\cite[Coro.~1.5]{Stoyanov11} (see also \cite{Liverani04} when $M$ is
compact, the result stated for the Liouville measure should extend to
the Bowen-Margulis measure).

\medskip The following result gives exponentially small error terms in
the equidistribution of the skinning measures to the Bowen-Margulis
measure, in the known situations when the geodesic flow is
exponentially mixing.  Here we state the result for skinning measures
but, clearly, it remains valid if $\sigma_\C$ is replaced by
$\sigma_{\flow{t}\Omega}$ as in Theorem \ref{theo:equid}.

\btheo\label{theo:exprate} Let $\wt M$ be a complete simply connected
Riemannian manifold with  negative sectional curvature. Let
$\Ga$ be a discrete, nonelementary group of isometries of $\wt
M$.  Let $\C=(C_i)_{i\in I}$ be a locally finite
$\Ga$-equivariant family of proper nonempty closed convex subsets of
$\wt M$, with finite nonzero skinning measure $\sigma_\C$.
Let $M=\Ga\backslash\wt M$.

\smallskip\noindent (i) If $M$ is compact and is $2$-dimensional or
locally symmetric, then there exist $\alpha\in\;]0,1[$ and
$\kappa''>0$ such that for all $\psi\in \operatorname{C}_{\rm
  c}^\alpha(T^1M)$, we have, as $t\ra+\infty$,
$$
\frac{1}{\|\sigma_{\C_t}\|}\int \psi \;d\sigma_{\C_t}=
\frac{1}{\|m_{\rm  BM}\|} \int \psi\;dm_{\rm BM}
+O(e^{-\kappa''t}\;\|\psi\|_\alpha)\;.
$$

\smallskip\noindent(ii) If $\wt M$ is a symmetric space and if $\Ga$
is an arithmetic lattice, then there exists $\ell\in\NN$ and
$\kappa''>0$ such that for all $\psi\in \C_{c}^\ell(T^1M)$, we have,
as $t\ra+\infty$,
$$
\frac{1}{\|\sigma_{\C_t}\|}\int \psi \;d\sigma_{\C_t}=
\frac{1}{\|m_{\rm  BM}\|} \int \psi\;dm_{\rm BM}
+O(e^{-\kappa''t}\;\|\psi\|_\ell)\;.
$$
\etheo

\dem Up to rescaling, we may assume that the sectional curvature is
bounded from above by $-1$. The critical exponent and the
Bowen-Margulis measure are finite in all cases considered.

Let us consider Claim (i).  Under these assumptions, there is some
$\alpha\in\;]0,1[\,$ such that the geodesic flow on $T^1M$ is
exponentially mixing for the H\"older regularity $\alpha$ and such
that the strong stable foliation of $T^1\wt M$ is $\alpha$-H\"older.

Fix $R>0$ and, for every $\eta>0$, let us consider the test function
$\phi_\eta=\phi_{\eta,\,R,\,\C}$ as in Section
\ref{sec:equidistribution}. Up to replacing $C_i$ by $\N_1C_i$, we may
assume that the boundary of $C_i$ is $C^{1,1}$-smooth, for every $i\in
I$ (see Section \ref{sec:geometry}).

Fix $\psi\in \operatorname{C}_{\rm c}^\alpha(T^1M)$. We may assume as
above that there exists a lift $\wt\psi: T^1\wt M\ra\RR$ of $\psi$
whose support is contained in a given fundamental domain $\Delta_\Ga$
for the action of $\Ga$ on $T^1\wt M$. First assume that $\Ga\bs I$ is
finite.  There exist $\eta_0>0$ and $t_0\geq 0$ such that for every
$\eta\in\;]0,\eta_0]$, and for every $t\in [t_0,+\infty[$, for every
$w\in T^1\wt M$ and $v\in V_{w,\,\eta,\,e^{-t}R}$, we have
\begin{equation}\label{eq:holdercont}
\wt\psi(v)= \wt\psi(w)+
\operatorname{O}\big((\eta+e^{-t})^\alpha\|\psi\|_\alpha\big)\;, 
\end{equation}
since
$d(v,w)=\operatorname{O}(\eta+e^{-t})$ by Equation
\eqref{eq:distlelongeod} and Lemma \ref{lem:compardisttunwss}.

As in the proof of Theorem \ref{theo:equid} using Equation
\eqref{eq:holdercont} instead of Equation \eqref{eq:unifcont} (see
Equation \eqref{eq:prepamixing}), we have
$$
\frac{\int \psi\;d\sigma_{\C_t}}{\|\sigma_{\C_t}\|}=\frac{\int_{T^1M}
\phi_\eta\;\psi\circ\flow t\;d\overline{m}_{\rm
  BM}}{\int_{T^1M}\phi_\eta \;d\overline{m}_{\rm BM}}
+\operatorname{O}\big((\eta+e^{-t})^\alpha\|\psi\|_\alpha\big)\;.
$$
As $M$ is compact, the Patterson densities and the Bowen-Margulis
measure are doubling measures and, using discrete convolution
approximation (see for instance \cite[p. 290-292]{Semmes96} or
\cite{KinKorShaTuo12}), there exists $\kappa'>0$ and, for every
$\eta>0$, a nonnegative function $\Phi_\eta \in \operatorname{C}_{\rm
  c}^\alpha(T^1M)$ such that
\begin{itemize}
\item $\int_{T^1M}\Phi_\eta\;d\overline{m}_{\rm BM}=\int_{T^1M}\phi_\eta
\;d\overline{m}_{\rm BM}$,
\item  $\int_{T^1M}|\Phi_\eta-\phi_\eta| \;d\overline{m}_{\rm BM}= 
\operatorname{O}(\eta\int_{T^1M}\phi_\eta
\;d\overline{m}_{\rm BM})$,
\item $\|\Phi_\eta\|_\alpha = \operatorname{O}(\eta^{-\kappa'}
\int_{T^1M}\phi_\eta \;d\overline{m}_{\rm BM})$.
\end{itemize}
Hence, applying the exponential mixing of the geodesic flow, with
$\kappa>0$ as in its definition, since $\int_{T^1M}\phi_\eta
\;dm_{\rm BM}$, which is equal to $\|\sigma_\C\|$ by Proposition
\ref{prop:integrable}, is independent of $\eta$, we have, for
$\eta\in\;]0,\eta_0]$ and $t\in [t_0,+\infty[$,
\begin{align*}
&\frac{\int \psi\;d\sigma_{\C_t}}{\|\sigma_{\C_t}\|}\\=&
\;\frac{\int_{T^1M}
\Phi_\eta\;\psi\circ\flow t\;d\overline{m}_{\rm BM}}{\int_{T^1M}\phi_\eta
\;d\overline{m}_{\rm BM}}
+ \operatorname{O}\big(\eta\; \|\psi\|_\infty+ 
(\eta+e^{-t})^\alpha\|\psi\|_\alpha\big)\\=\; & 
\frac{\int_{T^1M}  \Phi_\eta\;d\overline{m}_{\rm BM}}{\int_{T^1M}\phi_\eta
  \;d\overline{m}_{\rm BM}}
\;\int_{T^1M}  \psi\;d\overline{m}_{\rm BM}
+ \operatorname{O}\big(e^{-\kappa t}\|\Phi_\eta\|_\alpha\|\psi\|_\alpha+ 
\eta\; \|\psi\|_\infty+(\eta+e^{-t})^\alpha\|\psi\|_\alpha\big)\\=\; & 
\int_{T^1M}  \psi\;d\overline{m}_{\rm BM}
+ \operatorname{O}\big((e^{-\kappa t}\eta^{-\kappa'}+ 
\eta+(\eta+e^{-t})^\alpha)\|\psi\|_\alpha\big)\;.
\end{align*}
Taking $\eta=e^{-t\lambda}$ for $\lambda$ small enough (for instance
$\lambda =\kappa/(2\kappa')$ ), the result follows (for instance with
$\kappa''= \min\{\kappa/2,\; \kappa/(2\kappa'), \;\alpha \min\{1,
\kappa/(2\kappa')\}\}$ ), when $\Ga\bs I$ is finite.  As the implied
constants do not depend on the family $\C$, the result holds in
general.

The proof of Claim (ii) is similar. In this case, the strong stable
foliation is smooth and the Bowen-Margulis measure coincides, up to a
scalar multiple, with the Liouville measure. Thus, we can use the
usual convolution approximation (see for instance \cite[\S
1.6]{Ziemer89}) to approximate the test function by smooth functions.
\cqfd

{\small \bibliography{../biblio} }

\begin{thebibliography}{Bow2}

\bibitem[Bab]{Babillot02b}
M.~Babillot.
\newblock {\it On the mixing property for hyperbolic systems}.
\newblock {Israel J. Math. {\bf 129} (2002) 61--76}.

\bibitem[Bowd]{Bowditch95}
B.~Bowditch.
\newblock {\it Geometrical finiteness with variable negative curvature}.
\newblock {Duke Math. J. {\bf 77} (1995) 229--274}.

\bibitem[Bowe]{Bowen72a}
R.~Bowen.
\newblock {\it Periodic orbits for hyperbolic flows}.
\newblock {Amer. J. Math. {\bf 94} (1972), 1--30}.

\bibitem[BH]{BriHae99}
M.~R. Bridson and A.~Haefliger.
\newblock {\it Metric spaces of non-positive curva\-tu\-re}.
\newblock {Grund. math. Wiss. {\bf 319}, Springer Verlag, 1999}.

\bibitem[Bri]{Brin95}
M.~Brin.
\newblock {\it Ergodicity of the geodesic flow}.
\newblock {Appendix in W.~Ballmann, {\it Lectures on spaces 
of nonpositive curvature}, DMV Seminar {\bf 25}, Birkh\"auser,
1995, 81--95}.

\bibitem[Clo]{Clozel03}
L.~Clozel.
\newblock {\it D\'emonstration de la conjecture $\tau$}.
\newblock {Invent. Math. {\bf 151} (2003) 297--328}.

\bibitem[Dal1]{Dalbo99}
F.~Dal'Bo.
\newblock {\it Remarques sur le spectre des longueurs d'une surface et
  comptage}.
\newblock {Bol. Soc. Bras. Math. {\bf 30} (1999) 199--221}.

\bibitem[Dal2]{Dalbo00}
F.~Dal'Bo.
\newblock {\it Topologie du feuilletage fortement stable}.
\newblock {Ann. Inst. Fourier {\bf 50} (2000) 981--993}.

\bibitem[DOP]{DalOtaPei00}
F.~Dal'Bo, J.-P. Otal, and M.~Peign\'e.
\newblock {\it S\'eries de Poincar\'e des groupes g\'eom\'etriquement finis}.
\newblock {Israel J. Math. {\bf 118} (2000) 109--124}.

\bibitem[Dol]{Dolgopyat98}
D.~Dolgopyat.
\newblock {\it On decay of correlation in Anosov flows}.
\newblock {Ann. of Math. {\bf 147} (1998) 357--390}.

\bibitem[EM]{EskMcMul93}
A.~Eskin and C.~McMullen.
\newblock {\it Mixing, counting, and equidistribution in Lie groups}.
\newblock {Duke Math. J. {\bf 71} (1993) 181--209}.

\bibitem[GLP]{GiuLivPol12}
P.~Giulietti, C.~Liverani, and M.~Pollicott.
\newblock {\it Anosov flows and dynamical Zeta Functions}.
\newblock {Preprint {\tt [arXiv:1203.0904]}}.

\bibitem[Ham]{Hamenstadt89}
U.~Hamenst{\"a}dt.
\newblock {\it A new description of the Bowen-Margulis measure}.
\newblock {Erg. Theo. Dyn. Sys. {\bf 9} (1989) 455--464}.

\bibitem[HP1]{HerPau97}
S.~Hersonsky and F.~Paulin.
\newblock {\it On the rigidity of discrete isometry groups of negatively curved
  spaces}.
\newblock {Comm. Math. Helv. {\bf 72} (1997) 349--388}.

\bibitem[HP2]{HerPau04}
S.~Hersonsky and F.~Paulin.
\newblock {\it Counting orbit points in coverings of negatively curved
  manifolds and Hausdorff dimension of cusp excursions}.
\newblock {Erg. Theo. Dyn. Sys. {\bf 24} (2004) 1--22}.

\bibitem[HP3]{HerPau10}
S.~Hersonsky and F.~Paulin.
\newblock {\it On the almost sure spiraling of geodesics in negatively curved
  manifolds}.
\newblock {J. Diff. Geom. {\bf 85} (2010) 271--314}.

\bibitem[Kim]{Kim13}
I.~Kim.
\newblock {\it Counting, mixing and equidistribution of horospheres in
  geometrically finite rank one locally symmetric manifolds}.
\newblock {Preprint {\tt [arXiv:1103.5003]}}.

\bibitem[KKST]{KinKorShaTuo12}
J.~Kinnunen, R.~Korte, N.~Shanmugalingam and H.~Tuominen.
\newblock {\it A characterization of Newtonian functions with zero 
boundary values}.
\newblock {Calc. Var. Part. Diff. Eq. {\bf 43} (2012) 507--528}.

\bibitem[KM1]{KleMar96}
D.~Kleinbock and G.~Margulis.
\newblock {\it Bounded orbits of nonquasiunipotent flows on homogeneous
  spaces}.
\newblock {Sinai's Moscow Seminar on Dynamical Systems, 141--172, Amer. Math.
  Soc. Transl. Ser. {\bf 171}, Amer. Math. Soc. 1996}.

\bibitem[KM2]{KleMar99}
D.~Kleinbock and G.~Margulis.
\newblock {\it Logarithm laws for flows on homogeneous spaces}.
\newblock {Invent. Math. {\bf 138} (1999) 451--494}.

\bibitem[KO]{KonOh11}
A.~Kontorovich and H.~Oh.
\newblock {\it Apollonian circle packings and closed horospheres on hyperbolic
  3-manifolds}.
\newblock {J. Amer. Math. Soc. {\bf 24} (2011) 603--648}.

\bibitem[Liv]{Liverani04}
C.~Liverani.
\newblock {\it On contact Anosov flows}.
\newblock {Ann. of Math. {\bf 159} (2004) 1275--1312}.

\bibitem[Mar1]{Margulis69}
G.~Margulis.
\newblock {\it Applications of ergodic theory for the investigation of
  manifolds of negative curvature}.
\newblock {Funct. Anal. Applic. {\bf 3} (1969) 335--336}.

\bibitem[Mar2]{Margulis70}
G.~Margulis.
\newblock {\it Certain measures that are connected with $U$-flows on compact
  manifolds}.
\newblock {Funct. Anal. Applic. {\bf 4} (1970) 55--67}.

\bibitem[Mar3]{Margulis04}
G.~Margulis.
\newblock {\it On some aspects of the theory of Anosov systems}.
\newblock {Mono. Math., Springer Verlag, 2004}.

\bibitem[Mark]{Marklof10}
J.~Marklof.
\newblock {\it The asymptotic distribution of Frobenius numbers}.
\newblock {Invent. Math. {\bf 181} (2010) 179--207}.

\bibitem[OS1]{OhShaInv}
H.~Oh and N.~Shah.
\newblock {\it The asymptotic distribution of circles in the orbits of Kleinian
  groups}.
\newblock {Invent. Math. {\bf 187} (2012) 1--35}

\bibitem[OS2]{OhShaCounting}
H.~Oh and N.~Shah.
\newblock {\it Equidistribution and counting for orbits of geometrically finite
  hyperbolic groups}.
\newblock {To appear in J. Amer. Math. Soc.}.

\bibitem[OP]{OtaPei04}
J.-P. Otal and M.~Peign\'e.
\newblock {\it Principe variationnel et groupes kleiniens}.
\newblock {Duke Math. J. {\bf 125} (2004) 15--44}.

\bibitem[PP1]{ParPau12JMD}
J.~Parkkonen and F.~Paulin.
\newblock {\it \'Equidistribution, comptage et approximation par irrationnels
  quadratiques}.
\newblock {J. Mod. Dyn. {\bf 6} (2012) 1--40}.

\bibitem[PP2]{ParPauRev}
J.~Parkkonen and F.~Paulin.
\newblock {\it Counting  arcs in negative curvature}.
\newblock {Preprint {\tt [arXiv: 1203.0175]}}.

\bibitem[PP3]{ParPau13b}
J.~Parkkonen and F.~Paulin.
\newblock {\it Counting common perpendicular arcs in negative curvature}.
\newblock {In preparation}.

\bibitem[Pau]{Paulin13b}
F.~Paulin.
\newblock {\it Regards crois\'es sur les s\'eries de Poincar\'e et leurs
  applications}.
\newblock {Notes d'expos\'e, GDR Platon 3341 CNRS, Univ. Neuch\^atel, 7-9
  f\'evrier 2011, see
  http://www.math.u-psud.fr/$\sim$paulin/preprints/liste\_preprints.html}.

\bibitem[PPS]{PauPolSha11}
F.~Paulin, M.~Pollicott, and B.~Schapira.
\newblock {\it Equilibrium states in negative curvature}.
\newblock {Book in preparation}.

\bibitem[Rob1]{Roblin02}
T.~Roblin.
\newblock {\it Sur la fonction orbitale des groupes discrets en courbure
  n\'egative}.
\newblock {Ann. Inst. Fourier {\bf 52} (2002) 145--151}.

\bibitem[Rob2]{Roblin03}
T.~Roblin.
\newblock {\it Ergodicit\'e et \'equidistribution en courbure n\'egative}.
\newblock {M\'emoire Soc. Math. France, {\bf 95} (2003)}.

\bibitem[Sch]{Schapira04b}
B.~Schapira.
\newblock {\it Lemme de l'ombre et non divergence des horosph\`eres d'une
  vari\'et\'e g\'eom\'etriquement finie}.
\newblock {Ann. Inst. Fourier (Grenoble) {\bf 54} (2004) 939--987}.

\bibitem[Sem]{Semmes96}
S.~Semmes.
\newblock {\it Finding curves on general spaces through quantitative
  topology, with applications to Sobolev and Poincar\'e inequalities}.
\newblock {Selecta Math. {\bf 2} (1996) 155--295}.

\bibitem[Sto]{Stoyanov11}
L.~Stoyanov.
\newblock {\it Spectra of Ruelle transfer operators for axiom A flows}.
\newblock {Nonlinearity {\bf 24} (2011) 1089--1120}.

\bibitem[SV]{StrVel95}
B.~Stratmann and S.~L. Velani.
\newblock {\it The Patterson measure for geometrically finite groups with
  parabolic elements, new and old}.
\newblock {Proc. London Math. Soc. {\bf 71} (1995) 197--220}.

\bibitem[Sul]{Sullivan84}
D.~Sullivan.
\newblock {\it Entropy, Hausdorff measures old and new, and the limit set of
  geometrically finite Kleinian groups}.
\newblock {Acta Math. {\bf 153} (1984) 259--277}.

\bibitem[Wal]{Walter76}
R.~Walter.
\newblock {\it Some analytical properties of geodesically convex sets}.
\newblock {Abh. Math. Sem. Univ. Hamburg {\bf 45} (1976) 263--282}.

\bibitem[Zie]{Ziemer89}
W.~P.~Ziemer.
\newblock {\it Weakly differentiable functions}.
\newblock {Grad. Texts Math. {\bf 120}, Springer Verlag, 1989}.

\end{thebibliography}

\bigskip
{\small\noindent \begin{tabular}{l} 
Department of Mathematics and Statistics, P.O. Box 35\\ 
40014 University of Jyv\"askyl\"a, FINLAND.\\
{\it e-mail: jouni.t.parkkonen@jyu.fi}
\end{tabular}
\medskip

\noindent \begin{tabular}{l}
D\'epartement de math\'ematique, UMR 8628 CNRS, B\^at.~425\\
Universit\'e Paris-Sud,
91405 ORSAY Cedex, FRANCE\\
{\it e-mail: frederic.paulin@math.u-psud.fr}
\end{tabular}

}

\end{document}